\documentclass[a4paper,11pt,twoside]{amsart}

\usepackage[latin1]{inputenc}
\usepackage{amsmath}
\usepackage{amsthm}
\usepackage{amssymb}
\usepackage[all]{xy}
\usepackage{hyperref}

\newtheorem{thm}{Theorem}[section]

\newtheorem{prop}[thm]{Proposition}

\newtheorem{lem}[thm]{Lemma}
\newtheorem{cor}[thm]{Corollary}

\numberwithin{equation}{section}

\theoremstyle{definition}
\newtheorem{definition}[thm]{Definition}
\newtheorem{remark}[thm]{Remark}
\newtheorem{ex}[thm]{Example}

\newcommand{\qqed}{\hspace*{\fill}$\Box$}

\newcommand{\Db}{{\rm D}^{\rm b}}
\newcommand{\Dp}{{\rm D}_{\rm perf}}

\newcommand{\Aut}{{\rm Aut}}

\newcommand{\Pic}{{\rm Pic}}

\newcommand{\res}{{\rm res}}
\newcommand{\rk}{{\rm rk}}
\newcommand{\coh}{{\cat{Coh}}}

\newcommand{\End}{{\rm End}}
\newcommand{\Hom}{{\rm Hom}}
\newcommand{\Spec}{{\rm Spec}}
\newcommand{\Spf}{{\rm Spf}}

\newcommand{\dual}{\raisebox{-.5ex}{\makebox[.45em]{{\LARGE$\check{}$}}}}

\newcommand{\epi}{\twoheadrightarrow}

\newcommand{\cat}[1]{\begin{bf}#1\end{bf}}
\newcommand{\Ext}{{\rm Ext}}

\newcommand{\KE}{{\rm Ker}}

\renewcommand{\ker}{{\rm Ker}}

\newcommand{\Mod}[1]{{\ko_{#1}\text{-}\cat{Mod}}}
\newcommand{\cal}{\mathcal}
\newcommand{\ka}{{\cal A}}

\newcommand{\kc}{{\cal C}}
\newcommand{\kd}{{\cal D}}
\newcommand{\ke}{{\cal E}}
\newcommand{\kf}{{\cal F}}
\newcommand{\kg}{{\cal G}}
\newcommand{\kh}{{\cal H}}
\newcommand{\ki}{{\cal I}}
\newcommand{\kk}{{\cal K}}

\newcommand{\kn}{{\cal N}}

\newcommand{\ko}{{\cal O}}

\newcommand{\kt}{{\cal T}}
\newcommand{\kx}{{\cal X}}
\newcommand{\ky}{{\cal Y}}
\newcommand{\kz}{{\cal Z}}
\newcommand{\NN}{\mathbb{N}}
\newcommand{\ZZ}{\mathbb{Z}}
\newcommand{\QQ}{\mathbb{Q}}
\newcommand{\RR}{\mathbb{R}}
\newcommand{\CC}{\mathbb{C}}

\newcommand{\PP}{\mathbb{P}}
\newcommand{\XX}{\mathbb{X}}
\newcommand{\xx}{{\bf x}}
\newcommand{\yy}{{\bf y}}

\newcommand{\ddual}{\dual\dual}

\newcommand{\HH}{\mathrm{H\!H}}
\newcommand{\HO}{\mathrm{H}\Omega}
\newcommand{\HT}{\mathrm{HT}}

\renewcommand{\to}{\xymatrix@1@=15pt{\ar[r]&}}
\newcommand{\lto}{\xymatrix@1@=15pt{&\ar[l]}}
\renewcommand{\leftarrow}{\xymatrix@1@=15pt{&\ar[l]}}
\renewcommand{\rightarrow}{\xymatrix@1@=15pt{\ar[r]&}}
\renewcommand{\mapsto}{\xymatrix@1@=15pt{\ar@{|->}[r]&}}
\renewcommand{\epi}{\xymatrix@1@=19pt{\ar@{->>}[r]&}}
\renewcommand{\twoheadrightarrow}{\xymatrix@1@=19pt{\ar@{->>}[r]&}}
\renewcommand{\hookrightarrow}{\xymatrix@W=2pt@1@=15pt{\ar@{^(->}[r]&}}
\newcommand{\congpf}{\xymatrix@1@=15pt{\ar[r]^-\sim&}}
\renewcommand{\cong}{\simeq}

\begin{document}
\title{Derived equivalences of K3 surfaces and orientation}

\author[D.\ Huybrechts, E.\ Macr\`i, and P.\ Stellari]{Daniel Huybrechts, Emanuele Macr\`i, and Paolo Stellari}

\address{D.H.: Mathematisches Institut,
Universit{\"a}t Bonn, Beringstr.\ 1, 53115 Bonn, Germany}
\email{huybrech@math.uni-bonn.de}

\address{E.M.: Department of Mathematics, University of Utah, 155 South 1400 East, Salt
Lake City, UT 84112-0090, USA} \email{macri@math.utah.edu}

\address{P.S.: Dipartimento di Matematica ``F. Enriques'',
Universit{\`a} degli Studi di Milano, Via Cesare Saldini 50, 20133
Milano, Italy} \email{paolo.stellari@unimi.it}

\begin{abstract}
Every Fourier--Mukai equivalence between the derived categories of
two K3 surfaces induces a Hodge isometry of their cohomologies
viewed as Hodge structures of weight two endowed with the Mukai
pairing. We prove that this Hodge isometry preserves the natural
orientation of the four positive directions. This leads to a
complete description of the action of the group of all
autoequivalences on cohomology very much like the classical
Torelli theorem for K3 surfaces determining all Hodge isometries
that are induced by automorphisms.
\end{abstract}

\keywords{K3 surfaces, derived categories, deformations}

\subjclass[2000]{18E30, 14J28}

\maketitle



\section{Introduction}

The second cohomology $H^2(X,\ZZ)$ of a K3 surface $X$ is an even
unimodular lattice of signature $(3,19)$ endowed with a natural
weight two Hodge structure. The inequality $(\alpha,\alpha)>0$
describes an open subset of the $20$-dimensional real vector space
$H^{1,1}(X)\cap H^2(X,\RR)$ with two connected components $\kc_X$
and $-\kc_X$. Here $\kc_X$ denotes the positive cone, i.e.\ the
connected component that contains the K\"ahler cone $\kk_X$ of all
K\"ahler classes on $X$.

Any automorphism $f:X\congpf X$  of the complex surface $X$
defines an isometry $$\xymatrix{f_*:H^2(X,\ZZ)\ar[r]^-\sim&
H^2(X,\ZZ)}$$ compatible with the weight two Hodge structure. In
particular, $f_*$ preserves the set $\kc_X\sqcup(-\kc_X)$. As the
image of a K\"ahler class is again a K\"ahler class,  one actually
has $f_*(\kc_X)=\kc_X$. In other words, $f_*$ respects the
connected components of the set of $(1,1)$-classes $\alpha$ with
$(\alpha,\alpha)>0$. If one wants to avoid the existence of
K\"ahler structures, the proof of this assertion  is a little more
delicate. However, applying his polynomial invariants, Donaldson
proved in \cite{Donaldson}  a much stronger result not appealing
to the complex or K\"ahler structure  of $X$ at all.

Before recalling his result, let us rephrase the above discussion
in terms of orientations of positive three-spaces. Consider any
three-dimensional subspace $F\subset H^2(X,\RR)$ on which  the
intersection pairing is positive definite. Then $F$ is called a
positive three-space. Using orthogonal projections, given
orientations on two positive three-spaces can be compared to each
other. So, if $\rho$ is an arbitrary isometry of $H^2(X,\RR)$ and
$F$ is a positive three-space, one can ask whether a given
orientation of $F$ coincides with the image of this orientation on
$\rho(F)$. If this is the case, then one says that $\rho$ is
\emph{orientation preserving}. Note that this does neither depend
on $F$ nor  on the chosen orientation of $F$. The fact that any
automorphism $f$ of the complex surface $X$ induces a Hodge
isometry with $f_*(\kc_X)=\kc_X$ is equivalent to saying that
$f_*$ is orientation preserving. More generally one has:

\medskip

\noindent{\bf Theorem 1. (Donaldson)} \emph{Let $f:X\congpf X$ be any
diffeomorphism. Then the induced isometry $f_*:H^2(X,\ZZ)\congpf
H^2(X,\ZZ)$ is orientation preserving.}

\medskip

This leads to a complete description of the image of the natural
representation $$\xymatrix{{\rm Diff}(X)\ar[r]& {\rm
O}(H^2(X,\ZZ))}$$ as the set of all orientation preserving
isometries of the lattice $H^2(X,\ZZ)$. That every orientation
preserving isometry can be lifted to a diffeomorphism relies on
the Global Torelli theorem (see \cite{Borcea}). The other
inclusion is the above result of Donaldson.

\smallskip

There are several reasons to pass from automorphisms or, more
generally, diffeomorphisms of a K3 surface $X$ to derived
autoequivalences. First of all, exact autoequivalences of the
bounded derived category $$\Db(X):=\Db(\coh(X))\cong\Db_{\rm
coh}(\Mod{X})$$ of coherent sheaves can be considered as natural
generaliza\-tions of automorphisms of the complex surface $X$, for
any automorphism clearly induces an autoequivalence of $\Db(X)$.
The second motivation comes from mirror symmetry, which suggests a
link between the group of auto\-equi\-valences of $\Db(X)$ and the
group of diffeomorphisms or rather symplectomorphisms of the
mirror dual K3 surface.

In order to study the derived category $\Db(X)$ and its
autoequivalences, one needs to introduce the Mukai lattice
$\widetilde H(X,\ZZ)$ which comes with a natural weight two Hodge
structure. The lattice $\widetilde H (X,\ZZ)$ is by definition the
full cohomology $H^*(X,\ZZ)$ endowed with a modification of the
intersection pairing (the Mukai pairing) obtained by introducing a
sign in the pairing of $H^0$ with $H^4$. The weight two Hodge
structure on $\widetilde H(X,\ZZ)$ is by definition orthogonal
with respect to the Mukai pairing and therefore determined by
setting $\widetilde H^{2,0}(X):=H^{2,0}(X)$.

In his seminal article \cite{Mu}, Mukai showed that to any exact
autoequivalence (of Fourier--Mukai type)
$\Phi:\Db(X)\congpf\Db(X)$ of the bounded derived category of a
projective K3 surface there is naturally associated an isomorphism
$$\xymatrix{\Phi^{H^*}:\widetilde H(X,\ZZ)\ar[r]^-\sim& \widetilde
H(X,\ZZ)}$$ which respects the Mukai pairing and the Hodge
structure, i.e.\ $\Phi^{H^*}$ is a Hodge isometry of the Mukai
lattice. Thus, the natural representation $\Aut(X)\to{\rm
O}(H^2(X,\ZZ))$ is generalized to a representation
$$\xymatrix{\Aut(\Db(X))\ar[r]& {\rm
O}(\widetilde H(X,\ZZ)).}$$

The lattice $\widetilde H(X,\ZZ)$ has signature $(4,20)$ and, in
analogy to the discussion above, one says that an isometry $\rho$
of $\widetilde H(X,\ZZ)$ is \emph{orientation preserving} if under
orthogonal projection a given orientation of a positive four-space
in $\widetilde H(X,\ZZ)$ coincides with the induced  one on its
image under $\rho$. Whether $\rho$ is orientation preserving does
neither depend on the positive four-space nor on the chosen
orientation of it. The main result of this paper is the proof of a
conjecture that has been formulated by Szendr\H{o}i in \cite{Sz}
as the mirror dual of Donaldson's Theorem 1.

\medskip

\noindent{\bf Theorem 2.} \emph{Let $\Phi:\Db(X)\congpf\Db(X)$ be
an exact autoequivalence of the bounded derived category of a
projective K3 surface $X$. Then the induced Hodge isometry
$\Phi^{H^*}:\widetilde H(X,\ZZ)\congpf\widetilde H(X,\ZZ)$ is
orientation preserving.}

\medskip

For most of the known equivalences this can be checked directly,
e.g.\ for spherical twists and tensor products with line bundles.
The case of equivalences given by the universal family of stable
sheaves is more complicated and was treated in \cite{HS}. The
proof of the general case, as presented in this article, is rather
involved.

Theorem 2 can also be formulated for derived equivalences between
two different projective K3 surfaces by using the natural
orientation of the four positive directions (see Section
\ref{sect:noniso}).
\smallskip

Based on results of Orlov \cite{Or1}, it was proved in
\cite{HOY,P} that any orientation preserving Hodge isometry
actually occurs in the image of the representation ${\rm
Aut}(\Db(X))\to  {\rm O}(\widetilde H(X,\ZZ))$. This can be
considered as the analogue of the fact alluded to above that any
orientation preserving isometry of $H^2(X,\ZZ)$  lifts to a
diffeomorphism or to the part of the Global Torelli theorem that
describes the automorphisms of a K3 surfaces in terms of Hodge
isometries of the second cohomology. Together with Theorem 1, it
now allows one to describe the image of the representation
$\Phi\mapsto\Phi^{H^*}$ as the group of all orientation preserving
Hodge isometries of the Mukai lattice $\widetilde H(X,\ZZ)$:

\medskip

\noindent{\bf Corollary 3.} \emph{For any algebraic K3 surface $X$
one has} $${\rm Im}\left({\rm Aut}(\Db(X))\to {\rm O}(\widetilde
H(X,\ZZ))\right)={\rm O}_+(\widetilde H(X,\ZZ)).$$

\medskip

The kernel of ${\rm Diff}(M)\to {\rm O}(H^2(X,\ZZ))$ is largely
unknown, e.g.\ we do not know whether it is connected. In the
derived setting we have at least a beautiful conjecture due to
Bridgeland which describes the kernel of the analogous
representation in the derived setting as the fundamental group of
an explicit period domain (see \cite{B}).

\medskip

The key idea of our approach is actually quite simple: Deforming
the Fourier--Mukai kernel of a given derived equivalence yields a
derived equivalence between generic K3 surfaces and those have
been dealt with  in \cite{HMS}. In particular, it is known that in
the generic case the action on cohomology is orientation
preserving. As the action on the lattice $\widetilde H(X,\ZZ)$
stays constant under deformation, this proves the assertion.

What makes this program complicated and interesting, is the
deformation theory that is  involved. First of all, one has to
make sure that the Fourier--Mukai kernel does deform sideways to
any order. This can be shown if one of the two Fourier--Mukai
partners is deformed along a twistor space, which itself depends
on a chosen Ricci-flat metric on the K3 surface, and the other is
deformed appropriately. The second problem, as usual in
deformation theory, is convergence of the deformation. This point
is quite delicate for at least  two reasons: The Fourier--Mukai
kernel is not just a coherent sheaf but a complex of coherent
sheaves and the deformation we consider is not algebraic. We
circumvent both problems by deforming only to the very general
fibre of a formal deformation, which is a rigid analytic variety.
(In fact, only the abelian and derived category of coherent
sheaves on the rigid analytic variety are used and never the
variety itself.) The price one pays for passing to the general
fibre of the formal deformation only and not to an actual
non-algebraic K3 surface is that the usual $\CC$-linear categories
are replaced by categories defined over the non-algebraically
closed field $\CC((t))$ of Laurent series.

\medskip

The original paper \cite{HMS1} combined the results of this
article and the more formal aspects now written up in \cite{HMS3}.
We hope that splitting \cite{HMS1} in two shorter articles will
make the structure of the discussion clearer and not lead to
confusion.

The plan of this paper is as follows: In Section
\ref{sect:K3general}, after defining the formal setting we will
work with, we show that, for a formal twistor deformation
associated to a very general K\"ahler class, the bounded derived
category of its general fibre has only one spherical object up to
shift (Proposition \ref{prop:sphgen}). Hence the results of
\cite{HMS} can be applied. This part is based on results in
\cite{HMS1} not covered here, which can now also be found in the
separate \cite{HMS3}. In order to study autoequivalences of the
bounded derived category of the general fibre, we construct a
special stability condition for which the sections of the formal
deformations yield the only stable semi-rigid objects (Proposition
\ref{prop:semirig}). As a consequence, we prove that up to shift
and spherical twist any autoequivalence of the general fibre sends
points to points (Proposition \ref{prop:bijratpoints}) and its
Fourier--Mukai kernel is a sheaf (Proposition
\ref{prop:kernelissheaf}).

Section \ref{sect:defoFM} deals with the
deformation theory of kernels of Fourier--Mukai equivalences. In order to control the obstructions,
one has to compare the Kodaira--Spencer classes of the two sides
of the Fourier--Mukai equivalence, which will be done using the
language of Hochschild (co)homology. In particular we show that, under suitable hypotheses on the deformation and on the Fourier--Mukai kernel, the kernel itself deforms.

In Section \ref{sect:Proof} we come back to derived equivalences
of K3 surfaces and their deformations. We will prove in two steps
that the first order obstruction and all the higher order
obstructions are trivial. For one of the K3 surfaces the
deformation will be given by the twistor space and for the other
it will be constructed recursively. The conclusion of the proof of Theorem 2 is in Section \ref{sect:gospecial}.


\section{The very general twistor fibre of a K3
surface}\label{sect:K3general}

In this section we study very special formal deformations of
smooth projective K3 surfaces. The aim is to prove that the
derived category of what will be called the general fibre of the
formal deformation behaves similarly to the derived category of a
generic non-projective K3 surface.


\subsection{Formal deformations}

Let $R:=\CC[[t]]$ be the ring of power series in $t$ with field of fractions
$K:=\CC((t))$, the field of all Laurent series. For any $n$, the surjection $R\twoheadrightarrow R_n:=\CC[t]/(t^{n+1})$ yields a closed embedding
$\Spec(R_n)\subset\Spec(R)$, the $n$-th infinitesimal
neighbourhood of $0\in\Spec(R)$. The increasing sequence of closed subschemes
$0=\Spec(R_0)\subset \Spec(R_1)\subset\ldots\subset
\Spec(R_n)\subset\ldots$ defines the formal scheme $\Spf(R)$.

A \emph{formal deformation} of a smooth projective variety $X$ is a smooth and proper formal $R$-scheme $\pi:\kx\to\Spf(R)$, where $\kx$ is given by an inductive system of schemes
$\pi_n:\kx_n\to\Spec(R_n)$, smooth and proper over $R_n$, and isomorphisms
$$\kx_{n+1}\times_{R_{n+1}}\Spec(R_n)\cong\kx_n$$ over
$R_n$ such that $\kx_0=X$. While the topological space underlying
the scheme $\kx$ is $X$, the structure sheaf of $\kx$ is
$\displaystyle{\ko_\kx=\varprojlim\ko_{\kx_n}}$. For the rest of this
paper the natural inclusions will be denoted as follows ($m<n$):
\begin{eqnarray*}&\iota_n : \kx_n \hookrightarrow
\kx\qquad{\rm and}\qquad\iota:=\iota_0:X\hookrightarrow \kx;\\
&i_{m,n}: \kx_m \hookrightarrow \kx_n,~ i_n := i_{n,
n+1}:\kx_n\hookrightarrow\kx_{n+1}, ~~~{\rm
and}~~~j_n=i_{0,n}:X\hookrightarrow \kx_n.
\end{eqnarray*}

\begin{ex}\label{ex:twistor}
    Examples of formal deformations of a smooth projective variety $X$ are obtained by looking at smooth and proper families
    $\XX\to D$ of (usually non-algebraic) complex manifolds over a one-dimensional disk $D$ with local parameter
    $t$ and special fibre $X=\XX_0$. The infinitesimal neighbourhoods
    $\kx_n:=\XX\times_D\Spec(R_n)$, considered as
    $R_n$-schemes, form an inductive system and thus
    give rise to a formal $R$-scheme $\pi:\kx\to\Spf(R)$. Thus,
    although the nearby fibres $\XX_t$ of $X=\XX_0$ could be
    non-algebraic, the construction leads to the algebraic object
    $\kx$.

    If $X$ is a K3 surface, examples of such families are provided by the
    \emph{twistor space}
    $\pi:\XX(\omega)\to\PP(\omega)$ associated to a
    K\"ahler class $\omega$ on $X$. The total space $\XX(\omega)$ is a
    compact complex threefold, which is never algebraic nor K\"ahler (see
    \cite[Rem.\ 25.2]{GHJ}), and the projection $\pi$ is smooth and
    holomorphic onto the base $\PP(\omega)$, which is non-canonically
    isomorphic to $\PP^1$. The fibres are the complex manifolds
    obtained by hyperk\"ahler rotating the original complex structure
    defining $X$ in the direction of the hyperk\"ahler metric
    determined by $\omega$. In particular, there is a distinguished
    point $0\in\PP(\omega)$ such that the fibre
    $\XX(\omega)_0:=\pi^{-1}(0)$ is our original K3 surface $X$.
    By construction, the image of the composition
    \begin{equation*}\label{eqn:KSnorm}
    \xymatrix{T_0\PP(\omega)\ar[r]& H^1(X,\kt_X)\ar[r]&
    H^1(X,\Omega_X^1)}
    \end{equation*}
    of the Kodaira--Spencer map and the contraction $v\mapsto
    v\lrcorner\sigma=\sigma(v,-)$, where $\sigma\in H^0(X,\Omega_X^2)$
    is any non-trivial holomorphic two-form, is spanned by the
    K\"ahler class $\omega$ (for further details, see \cite{Beau}).

    Choosing a local parameter $t$ around $0$, one gets a
    formal deformation $\pi:\kx\to\Spf(R)$ which we call the \emph{formal twistor space} of $X$.
    Notice that the construction depends on the choice of the K\"ahler class
    $\omega$ and of the local parameter $t$.
\end{ex}

The $R$-linear category $\coh(\kx)$ of coherent
sheaves on $\kx$ contains the full abelian subcategory $\coh(\kx)_0\subset \coh(\kx)$
consisting of all sheaves $E\in\coh(\kx)$ such that $t^n E=0$ for $n\gg0$. (For the definition of coherent sheaves on noetherian formal
schemes see \cite[Ch.\ II.9]{HartAG} or \cite{IllFGA}.)
By definition $\coh(\kx)_0$
is a Serre subcategory and the quotient category
$$\coh(\kx_K):=\coh(\kx)/\coh(\kx)_0$$ is called
the \emph{category of coherent sheaves on the general fibre}. By
abuse of notation, we sometimes denote $\coh(\kx_K)$ by $\kx_K$.
When $\kx$ is a formal twistor space, we call the general fibre
$\kx_K$ the \emph{general twistor fibre}.

For $E\in\coh(\kx)$, denote by $E_K$ its projection in $\coh(\kx_K)$. The category $\coh(\kx_K)$ is a $K$-linear abelian category and
\begin{equation}\label{eq:morcoh}
\Hom_{\kx_K}(E_K,F_K)\cong\Hom_\kx(E,F)\otimes_R K,
\end{equation}
for any $E,F\in\coh(\kx)$ (see \cite[Prop.\ 2.4]{HMS1} or
\cite[Prop.\ 2.3]{HMS3}).

A coherent sheaf $E\in\coh(\kx)$ is \emph{$R$-flat} if the
multiplication with $t$ yields an injective homomorphism $t:E\to
E$. We denote by $\coh(\kx)_{\rm f}\subset\coh(\kx)$ the full
additive subcategory of all $R$-flat sheaves. Observe that
$\coh(\kx)_{\rm f}\to\coh(\kx_K)$ is essentially surjective, i.e.\
every object $F\in\coh(\kx_K)$ can be lifted to an $R$-flat sheaf
on $\kx$. Indeed, if $F=E_K$, then $(E_{\rm f})_K\cong E_K=F$,
where $E_\mathrm{f}=E/T$ and $T:=\cup\KE(t^n:E\to E)$ (notice that
since locally a coherent sheaf $E$ is the completion
of a finitely generated module over a noetherian ring, the union
stabilizes). By definition $E_{\rm
f}$ is an $R$-flat lift of $F$.

\bigskip

Passing to derived categories, consider the full thick triangulated subcategory
\[
\Db_0(\kx)\subseteq\Db(\kx):=\Db_{\rm coh}(\Mod{\kx})
\]
consisting of complexes of $\ko_\kx$-modules with cohomologies in $\coh(\kx)_0$. The Verdier quotient category $\Db(\kx_K):=\Db(\kx)/\Db_0(\kx)$ is called the \emph{derived category of the general fibre}. As before, we denote by $E_K$ the projection to $\Db(\kx_K)$ of any $E\in\Db(\kx)$. The category $\Db(\kx_K)$ is a $K$-linear triangulated category and
\begin{equation*}\label{eq:morder}
\Hom_{\kx_K}(E_K,F_K)\cong\Hom_\kx(E,F)\otimes_R K,
\end{equation*}
for any $E,F\in\Db(\kx)$ (see \cite[Prop.\ 3.9]{HMS1} or
\cite[Prop.\ 2.9]{HMS3}). In particular $\Db(\kx_K)$ has finite
dimensional Hom-spaces over $K$. Moreover $\coh(\kx_K)$ is the
heart of a bounded $t$-structure in $\Db(\kx_K)$.

When $X$ is a K3 surface, the main properties of $\Db(\kx_K)$ are
summarized by the following result which is proved in \cite[Sect.\
3]{HMS1} or \cite[Thm.\ 1.1]{HMS3}.

\begin{thm}\label{prop:K3cat} Let $\pi:\kx\to\Spf(R)$ be a formal deformation of a K3
surface $X=\kx_0$. Then the derived category $\Db(\kx_K)$ of its
general fibre is equivalent to $\Db(\coh(\kx_K))$ which is a $K$-linear K3 category.
\end{thm}

Recall that a \emph{K3 category} is a triangulated category with finite dimensional Hom-spaces and such that the double shift defines a Serre functor (see \cite{HMS}).

\bigskip

The main derived functors (tensor product, pull-back push-forward,
Hom's) are well-defined at the level of derived categories of
formal schemes over $\Spf(R)$. Moreover, all the basic properties
of them (e.g.\ commutativity, flat base change, projection
formula) hold in the formal context. All those functors are
$R$-linear and hence they factorize through the derived category
of the general fibre, verifying the same compatibilities (see
\cite[App.\ A.1]{HMS1} or \cite[Sect.\ 2.3]{HMS3}). To simplify
the notation, sometimes we will denote a functor and its derived
version in the same manner. In the case of an immersion of
(formal) schemes $j:Y\hookrightarrow Z$ and a sheaf $F\in\coh(Z)$,
we set $F|_Y:=\kh^0(Lj^*F)$.

\begin{remark}\label{rmk:functors}
    Let $X$ and $Y$ be smooth and projective varieties. Let $\kx_n,\ky_n\to\Spec(R_n)$ be an inductive system of smooth and proper schemes such that $\kx_n\times_{R_n}\Spec(R_0)\cong X$ and $\ky_n\times_{R_n}\Spec(R_0)\cong Y$, with $n\in\NN$.
    These collections yield  formal deformations $\kx,\ky\to\Spf(R)$ of $X$ and $Y$ respectively.

    i) Any bounded complex with coherent cohomology on a smooth
    formal scheme is \emph{perfect}, i.e.\ locally quasi-isomorphic to a
    finite complex of locally free sheaves of finite type.
    This is however not true for $\kx_n$, $n>0$, in which case we will have sometimes to work with $$\Dp(\kx_n)\subset\Db(\kx_n),$$ the full triangulated subcategory of perfect complexes on $\kx_n$. For $E\in\Db(\kx)$, we set $E_n:=L\iota_n^* E\in\Dp(\kx_n)$.

    ii) Given $\ke\in\Db(\kx\times_R\ky)$, we can define the \emph{Fourier--Mukai transform}
    $$\xymatrix{\Phi_\ke:\Db(\kx)\ar[r]&\Db(\ky)},~~~\xymatrix{E\ar@{|->}[r]&Rp_*(q^*E\otimes^L\ke),}$$
    where $p:\kx\times\ky\to\ky$ and $q:\kx\times\ky\to\kx$ are the projections. All the basic properties of Fourier--Mukai transforms valid for smooth
    projective varieties extend to the formal setting (see \cite[App.\ A.2]{HMS1} or \cite[Sect.\ 2.3]{HMS3}).
    Everything said also works for
    $\kx_n$ and $\ky_n$ with the only difference that we have to assume now that the Fourier--Mukai
    kernel $\ke_n\in\Db(\kx_n\times_{R_n}\ky_n)$ is perfect. Analogously, one can define the Fourier--Mukai transform
    $$\xymatrix{\Phi_\kf:\Db(\kx_K)\ar[r]&\Db(\ky_K)}$$
    associated to an object $\kf\in\Db((\kx\times_R\ky)_K)$. Indeed, given
    $\ke\in\Db(\kx\times_R\ky)$ with $\ke_K\cong\kf$, by $R$-linearity, the Fourier--Mukai transform
    $\Phi_\ke:\Db(\kx)\to\Db(\ky)$ descends to a Fourier--Mukai transform $\Phi_\kf:\Db(\kx_K)\to\Db(\ky_K)$, i.e.\ one has a
    commutative diagram $$\xymatrix{\Db(\kx)\ar[d]\ar[r]^{\Phi_\ke}&\Db(\ky)\ar[d]\\\Db(\kx_K)\ar[r]_{\Phi_\kf}&\Db(\ky_K).}$$

    iii) Let $\ke_n\in\Dp(\kx_n\times_{R_n}\ky_n)$, with
    $n\in\NN$, be such that its restriction $\ke_0:=Lj^*_n\ke_n\in\Db(X\times Y)$ is the kernel of a Fourier--Mukai
    equivalence $\Phi_{\ke_0}:\Db(X)\congpf\Db(Y)$. Then the
    Fourier--Mukai transforms $\Phi_{\ke_n}:\Dp(\kx_n)\to\Dp(\ky_n)$ and $\Phi_{\ke_n}:\Db(\kx_n)\to\Db(\ky_n)$
    are equivalences. The same holds true for $\ke\in\Db(\kx\times_R\ky)$ (see \cite[Prop.\ 3.19]{HMS1} or \cite[Prop.\ 2.12]{HMS3}).

    iv) As a consequence of iii), for $\ke\in\Db(\kx\times_R\ky)$ such that $\Phi_{\ke_0}:\Db(X)\congpf\Db(Y)$ is an equivalence, the
    Fourier--Mukai transform $\Phi_{\ke_K}:\Db(\kx_K)\congpf\Db(\ky_K)$ is an equivalence, where $\ke_K\in\Db((\kx\times_R\ky)_K)$ (see \cite[Cor.\ 3.20]{HMS1} or \cite[Cor.\ 2.13]{HMS3}).
\end{remark}

Sheaves and complexes of sheaves are usually denoted by $E$, $F$, etc. The use of $\ke$, $\kf$ wants to indicate that they are Fourier--Mukai kernels, which we wish to distinguish from the objects on the source and target variety of the associated Fourier--Mukai transform.

\bigskip

Given a formal deformation $\pi:\kx\to\Spf(R)$ of a smooth projective variety $X$, the category $\Db(\kx_K)$ contains the special object $\ko_{\kx_K}:=(\ko_\kx)_K$. Other objects of interest for this paper are obtained as follows. A
\emph{multisection} is an integral formal subscheme
$\kz\subset \kx$ which is flat of relative dimension zero over
$\Spf(R)$. The structure sheaf $\ko_\kz$ of such a multisection
induces an object in $\coh(\kx_K)$. Objects of this form will
usually be denoted by $K(\xx)\in\coh(\kx_K)$ and should be thought
of as (structure sheaves of) closed points $\xx\in\kx_K$ of the general fibre $\kx_K$. By
specialization, any point $K(\xx)\in\coh(\kx_K)$ determines a
closed point $x\in X$ of the special fibre.
The point $\xx$ is called $K$-\emph{rational} if $\kz\subset\kx$
is a section, i.e.\  $\pi|_\kz:\kz\to\Spf(R)$ is an isomorphism. Clearly, a closed point $\xx\in\kx_K$ is $K$-rational if and only if the natural
homomorphism $K\to\End_{\kx_K}(K(\xx))$ is an isomorphism.

\begin{remark}\label{rmk:ratpoints} Let $\kx$ be a formal deformation of a K3 surface $X$ and take $F\in\Db(\kx_K)$.

i) We call $F$ \emph{rigid} if $\Ext_{\kx_K}^1(F,F)=0$. We call $F$ \emph{spherical} if it is rigid and $\Ext_{\kx_K}^i(F,F)\cong K$ if $i=0,2$. An example of spherical object is provided by $\ko_{\kx_K}$. For this use \eqref{eq:semicont} and \eqref{eq:chi}.

ii) The object $F$ is \emph{semi-rigid} if $\Ext_{\kx_K}^1(F,F)=K^{\oplus 2}$. Applying again \eqref{eq:chi} and Serre duality one shows that a $K$-rational point $\xx$ is semi-rigid.
\end{remark}


\subsection{Torsion (free) sheaves on the general fibre}\label{Sect:Tor}

For a formal deformation $\kx$ of a smooth projective variety $X$, we say that $F\in\coh(\kx_K)$ is
\emph{torsion} (resp.\ \emph{torsion free}) if there exists a lift
$E\in\coh(\kx)$ of $F$ which is a torsion (resp.\ {torsion free})
sheaf on $\kx$.

Note that $F\in\coh(\kx_K)$ is torsion if and only if any lift of
$F$ is torsion. A torsion free $F$ always admits also lifts which
are not torsion free (just add $R$-torsion sheaves). However, the
lift $E$ of a torsion free $F$ is $R$-flat if and only if it is
torsion free. We leave it to the reader to show that any subobject
of a torsion free $F\in\coh(\kx_K)$ is again torsion free and that
any $F\in\coh(\kx_K)$ admits a maximal torsion subobject $F_\mathrm{tor}\subset F$ whose cokernel $F/F_\mathrm{tor}$ is torsion free
(use \eqref{eq:morcoh}).

Well-known arguments of Langton  and Maruyama  can be adapted to
prove the following:

\begin{lem}\label{lem-Langton} Any torsion free $F\in\coh(\kx_K)$ admits an
$R$-flat lift $E\in\coh(\kx)$ such that the restriction $E_0$ of
$E$ to the special fibre is a torsion free sheaf on $X$.
\end{lem}

\begin{proof} We shall prove the following more precise
claim (cf.\ the proof of \cite[Thm.\ 2.B.1]{HL}): Let $E$ be a
torsion free (as $\ko_\kx$-module) coherent sheaf on $\kx$. Then
there exists a coherent subsheaf $E'\subset E$ with $E'_0:=L\iota^* E'$ a torsion
free sheaf and such that the inclusion induces an isomorphism $E'_K\cong
E_K$.

Suppose there is no such $E'\subset E$. Then we construct a
strictly decreasing sequence $ \ldots E^{n+1}\subset
E^n\subset\ldots \subset E^0=E$ inductively as follows:
$$E^{n+1}:=\ker\left(E^n\to E^n_0\to E_0^n/(E_0^n)_{\rm tor}\right),$$
where $(E_0^n)_{\rm tor}$ means the torsion part on the
special fibre. Clearly, $E_K^n=E_K$. For later use, we introduce
$B^n:=(E_0^n)_{\rm tor}$ and $G^n:=E_0^n/B^n$, which will be
considered simultaneously as sheaves on the special fibre $X$ and
as sheaves on $\kx$.
 Then there are two
exact sequences of sheaves on $\kx$, respectively $X$
$$\xymatrix{0\ar[r]&E^{n+1}\ar[r]&E^n\ar[r]&G^n\ar[r]&0}$$and
$$\xymatrix{0\ar[r]&G^{n}\ar[r]&E^{n+1}_0\ar[r]&B^n\ar[r]&0.}$$
The first exact sequence is just the definition of $E^{n+1}$. For
the second one we first construct the surjection $E_0^{n+1}\epi
B^n$ by restricting $E^{n+1}\to E^n$ to the closed fibre. The
image of the resulting homomorphism $E_0^{n+1}\to E_0^n$ is the
kernel of $E_0^n\to G^n$, i.e.\ $B^n$. Let $K$ be the kernel of
$E_0^{n+1}\to B^n$. We will show $K\cong G^n$. For this we use the
two short exact sequences
$$\xymatrix{0\ar[r]&tE^{n+1}\ar[r]&E^{n+1}\ar[r]&E_0^{n+1}\ar[r]&0}$$and
$$\xymatrix{0\ar[r]&tE^{n}\ar[r]&E^{n+1}\ar[r]&B^n\ar[r]&0,}$$
where the latter one is obtained from snake lemma applied to the
natural surjection $E^n/tE^n=E_0^n\epi E^n/E^{n+1}=G^n$. Another
diagram chase  shows that $K$ sits inside the short exact sequence
\[
\xymatrix{0\ar[r]&tE^{n+1}\ar[r]^{\varphi}&tE^{n}\ar[r]&K\ar[r]&0,}
\]
where $\varphi$ is the morphism induced by the inclusion
$E^{n+1}\to E^n$ via the isomorphisms $E_n\cong tE_n$ and $E_{n+1}\cong tE_{n+1}$, given by the multiplication by $t$. Thus $K\cong G^n$.

As $G^n$ is a torsion free sheaf on $X$, one has $B^{n+1}\cap
G^n=0$ in $E^{n+1}_0$. Therefore, there is a descending filtration of torsion
sheaves $\ldots \subset B^{n+1}\subset B^n\subset \ldots$ and an
ascending sequence of torsion free sheaves $\ldots \subset
G^n\subset G^{n+1}\subset\ldots$. The  support of the torsion
sheaves $B^n$ on $X$ might have components of codimension one, but
for $n\gg0$ the filtration stabilizes in codimension one. Indeed, clearly, the support of the $B^n$ stabilizes for
$n\gg0$ and on there the generic rank will have to stabilize,
which then means that the $B^n$ themselves stabilize in
codimension one.

Hence,
$\ldots\subset G^n\subset G^{n+1}\subset\ldots$ stabilizes for
$n\gg0$ in codimension one as well. In particular, the reflexive
hulls do not change, i.e.\ $(G^n)\ddual=(G^{n+1})\ddual$ for
$n\gg0$. Therefore, for $n\gg0$ the sequence $G^n\subset
G^{n+1}\subset\ldots$ is an ascending sequence of coherent
subsheaves of a fixed coherent sheaf and hence stabilizes  for
$n\gg0$. This in turn implies that $\ldots \subset B^{n+1}\subset
B^n\subset\ldots$ stabilizes for $n\gg0$.

Replacing $E$ by $E^n$ with $n\gg0$, we may assume that
$G:=G^0=G^1=\ldots=G^n=\ldots$ and $0\ne
B:=B^0=B^1=\ldots=B^n=\ldots$. Note that this actually implies
$E_0=G\oplus B$.

We continue with the new $E$ obtained in this
way and consider the filtration $E^n$ for it. Now set
$Q^n:=E/E^n$. Then by definition of $E^n$ one has
$Q_0^n\cong G$. Moreover, there exists an exact sequence
$$\xymatrix{0\ar[r]&G\ar[r]&Q^{n+1}\ar[r]&Q^n\ar[r]&0,}$$
for $E^n/E^{n+1}\cong G^n=G$. Next, the quotient
$E\twoheadrightarrow Q^n$ factorizes over $E\twoheadrightarrow
E/t^nE\twoheadrightarrow Q^n$. Indeed, by construction
$tE^n\subset E^{n+1}$ and thus $t^nE=t^nE^0\subset E^n$. Thus, we
have a sequence of surjections $E/t^nE\twoheadrightarrow Q^n$ of
coherent sheaves on $\kx_{n-1}$ whose restriction to the special
fibre yields the surjection $E_0\twoheadrightarrow G$ with
non-trivial torsion kernel $B$.

One easily verifies that the system $(E/t^nE\twoheadrightarrow
Q^n)$ yields a surjection $E\twoheadrightarrow Q$ of coherent
sheaves on the formal scheme $\kx$. Indeed, the system $(Q^n)$
defines a coherent sheaf on the formal scheme $\kx$, for
$G=\ker(Q^{n+1}\to Q^n)=t^n Q^{n+1}$. The inclusion
$t^nQ^{n+1}\subset G$ is obvious and $G\subset t^nQ^{n+1}$  can be
proved inductively as follows: Suppose one has proved already that
$G\subset t^kQ^{n+1}$ for $k<n$, is the kernel of the projection
$t^kQ^{n+1}\to t^kQ^n$. Then use $t^kQ^{n+1}/t^{k+1}Q^{n+1}\cong
t^kQ^n/t^{k+1}Q^n\cong G$ to deduce that $G=\ker(t^{k+1}Q^{n+1}\to
t^{k+1}Q^n)$. The compatibility with the quotient maps
$E/t^nE\twoheadrightarrow Q^n$ is obvious.

Outside the support of $B$ the morphism $E\twoheadrightarrow Q$ is
an isomorphism and hence $\ker(E\twoheadrightarrow Q)$ must be
torsion and non-trivial. This contradicts the assumption on $E$.
\end{proof}


\subsection{The $K$-group of the general fibre}\label{subsec:genNS}

With the usual notation, for $E,E'\in\Db(X)$ one sets:
$$\chi_0(E,E'):=\sum(-1)^i\dim_\CC\Ext_X^i(E,E')$$
and analogously for $F,F'\in\Db(\kx_K)$:
$$\chi_K(F,F'):=\sum(-1)^i\dim_K\Ext_{\kx_K}^i(F,F').$$
Recall that (see \cite[Cor.\ 3.15 and 3.16]{HMS1} or \cite[Cor.\
3.2 and 3.3]{HMS3}) for $E,E'\in\Db(\kx)$, one has
\begin{equation}\label{eq:chi}
    \chi_0(E_0,E'_0)=\chi(E_K,E'_K)
\end{equation}
and the following semi-continuity result
\begin{equation}\label{eq:semicont}
    \dim_\CC\Hom_X(E_0,E'_0)\geq\dim_K\Hom_{\kx_K}(E_K,E'_K),
\end{equation}
where $E_0$ and $E'_0$ are the restrictions of $E$ and $E'$ to the special fibre.

\medskip

Let us now consider the $K$-groups of the various derived
categories: $K(X):=K(\Db(X))$, $K(\kx):=K(\Db(\kx))$, and
$K(\kx_K):=K(\Db(\kx_K))$. We say that a class $[E]\in K(X)$ is \emph{numerically
trivial}, $[E]\sim 0$, if $\chi_0(E,E')=0$ for all $E'\in\Db(X)$.
Numerical equivalence for the general fibre is defined similarly
in terms of $\chi_K$. Set
$$\kn(X):=K(X)/\sim \qquad {\rm and} \qquad \kn(\kx_K):=K(\kx_K)/\sim.$$

\begin{lem}\label{lem:res}
Sending $[F]=[E_K]\in \kn(\kx_K)$ (where $E\in\Db(\kx)$ is any lift of $F$) to $[E_0]\in\kn(X)$ determines
an injective linear map
$$\xymatrix{\res:\kn(\kx_K)\ar[r]&\kn(X)/\iota^*K(\kx)^\perp.}$$
(The orthogonal complement is taken with respect to $\chi_0$.)
\end{lem}

\begin{proof} The linearity of the map is evident, but in order to show that it  is well-defined one needs that  $\chi_K(E_K,-)\equiv 0$ implies
$\chi_0(E_0,E'_0)=0$ for all $E'\in\Db(\kx)$. This follows from \eqref{eq:chi}.

In order to prove injectivity of $\res$, suppose
$E_0\in\iota^*K(\kx)^\perp$. Then $\chi_K(E_K,E'_K)=\chi_0(E_0,E'_0)=0$
for all $[E']\in K(\kx)$. Since $K(\kx)\to K(\kx_K)$ is
surjective, this proves the claim.
\end{proof}

\begin{remark}\label{rem:resonK}
In fact, $\res$ can be lifted to a map
$$\xymatrix{K(\kx_K)\ar[r]& K(X),}$$
which will be used only once (see the proof of Corollary
\ref{cor:sheafcannot}). To show that the natural map $[E_K]\mapsto
[E_0]$ is well-defined, it suffices to show that any $R$-torsion
sheaf $E\in\coh(\kx)$ leads to a trivial class $[L\iota^*E]=[E_0]$
in $K(X)$.

As any $R$-torsion sheaf admits a filtration with quotients living
on $\kx_0=X$, it is enough to prove that $0=[L\iota^*\iota_*G]\in
K(X)$ for any $G\in \coh(X)$. For this, we complete the adjunction
morphism $L\iota^*\iota_*G\to G$ to the distinguished triangle
\begin{equation}\label{disp:ses}
\xymatrix{G[1]\ar[r]&L\iota^*\iota_* G\ar[r]&G},
\end{equation}
which shows $[L\iota^*\iota_*G]=[G]+[G[1]]=0$. For the existence
of \eqref{disp:ses} see e.g.\ \cite[Cor.\ 11.4]{HFM}. The proof
there can be adapted to the formal setting.
\end{remark}


\subsection{The general fibre of a very general twistor space}\label{sect:Twistor}

Let $\pi:\kx\to\Spf(R)$ be a formal twistor space of a K3 surface
$X$ associated to a K\"ahler class $\omega$. In the following,
$\omega$ has to be chosen
very general in order to ensure that only
the trivial line bundle $\ko_X$ deforms sideways. Here is the
precise definition we shall work with.

\begin{definition}
A K\"ahler class $\omega\in H^{1,1}(X,\RR)$ is called \emph{very
general} if there is no non-trivial integral class $0\ne\alpha\in
H^{1,1}(X,\ZZ)$ orthogonal to $\omega$, i.e.\ $\omega^\perp\cap
H^{1,1}(X,\ZZ)=0$.
\end{definition}

The twistor space associated to a very general K\"ahler class will
be called a \emph{very general twistor space} and its general
fibre  a \emph{very general twistor fibre}.

\begin{remark} Thus the set of very general K\"ahler classes is the complement
(inside the K\"ahler cone $\kk_X$) of the countable union of all
hyperplanes $\alpha^\perp\subset H^{1,1}(X,\RR)$ with
$0\ne\alpha\in H^{1,1}(X,\ZZ)$ and is, therefore, not empty. Moreover,
very general K\"ahler classes always exist also in
$\Pic(X)\otimes\RR$.
\end{remark}

In the next proposition we collect the consequences of this choice
that will be used in the following discussion.

\begin{prop}\label{prop:omegagen} Let $\omega$ be a very general K\"ahler class on a K3 surface $X$ and
let $\pi:\kx\to\Spf(R)$ be the induced formal twistor space.
Then the following conditions hold:

{\rm i)} Any line bundle on $\kx$ is trivial.

{\rm ii)} Let $\kz\subset\kx$ be an $R$-flat formal subscheme.
Then either  the projection $\pi:\kz\to\Spf(R)$ is of relative
dimension zero or $\kz=\kx$.

{\rm iii)} The Mukai vector $v:=\mathrm{ch}\cdot\sqrt{\mathrm{td}(X)}$ and the restriction map $\res$
(see Lemma \ref{lem:res}) define  isomorphisms
$$\xymatrix{v\circ \res:\kn(\kx_K)\ar[r]^-\sim&\kn(X)/\iota^*K(\kx)^\perp\ar[r]^-\sim&(H^0\oplus H^4)(X,\ZZ)\cong\ZZ^{\oplus 2}.}$$
\end{prop}

\begin{proof}
i) As $\omega$ is very general, even to first order no integral
$(1,1)$-class on $X$ stays pure. Thus, in fact any line bundle  on
$\kx_1$ is trivial (see e.g.\ \cite[Lemma 26.4]{GHJ}).

ii) The second assertion holds without any genericity assumption
on the K\"ahler class $\omega$ and goes back to Fujiki
\cite{Fujiki}. The case of relative dimension one can also be
excluded using i).

iii)  For the second isomorphism we need to show that
$\iota^*K(\kx)^\perp/\!\!\sim\,=H^{1,1}(X,\ZZ)$. As has been used
already in the proof of i), no class in $H^{1,1}(X,\ZZ)$ deforms
even to first order. In other words, the image of
$\iota^*K(\kx)$ in $\kn(X)$ is contained in $H^0(X,\ZZ)\oplus
H^4(X,\ZZ)$ and thus orthogonal to $H^2(X)$ and in particular to
$H^{1,1}(X,\ZZ)$, proving
$H^{1,1}(X,\ZZ)\subset\iota^*K(\kx)^\perp$.

To prove that the inclusion
$H^{1,1}(X,\ZZ)\subset\iota^*K(\kx)^\perp$ is an equality,
consider $\ko_\kx$ and the structure sheaf $\ko_{L_x}$ of any
section through a given closed point $x\in X$. Then
$\iota^*\ko_\kx\cong\ko_X$ and $\iota^*\ko_{L_x}\cong k(x)$ with
Mukai vectors $(1,0,1)$ and $(0,0,1)$, respectively. These two
vectors form a basis of $(H^0\oplus H^4)(X,\ZZ)$ and their images
in $\kn(X)/\iota^*K(\kx)^\perp$ are linearly independent, because
$\chi_0(k(x),\iota^*\ko_{L_x})=\chi_0(k(x),k(x))=0$ but
$\chi_0(\ko_X,\iota^*\ko_{L_x})\ne0$. This proves that the
inclusion $H^{1,1}(X,\ZZ)\subset\iota^*K(\kx)^\perp$ cannot be
strict. Hence we get the second isomorphism.

The injectivity of the map $\res$ has been shown in general in
Lemma \ref{lem:res} and $[\ko_X]$ and $[k(x)]$, spanning
$\kn(X)/\iota^*K(\kx)^\perp$, are clearly in the image of it.
\end{proof}

\begin{ex}\label{exas:MukaiVect} Under the assumptions of
Proposition \ref{prop:omegagen},  we often write $(r,s)$ instead
of $(r,0,s)$ for the Mukai vector in the image of $v\circ\res$.

i) If $F$ is a non-trivial torsion free sheaf on $\kx_K$, then
$v(\res(F))=(r,s)$ with $r>0$.

ii) For any closed point $\yy\in\kx_K$ one has
$v(\res(K(\yy)))=(0,d)$, where $d$ is the degree (over $\Spf(R)$)
of the multisection $\kz\subset\kx$ corresponding to $\yy$.

iii) If $F\in\coh(\kx_K)$ with $v(\res(F))=(0,0)$, then $F=0$.
Indeed, if $E$ is an $R$-flat lift of $F$, then $E_0$ would be a
sheaf concentrated in dimension zero without global sections.
Hence $E_0=0$ and then also $E=0$.
\end{ex}

The restriction  $E_0$ of an $R$-flat lift $E$ of a torsion
$F\in\coh(\kx_K)$ is a torsion sheaf on the special fibre $\kx_0\cong X$
with zero-dimensional support (use Proposition
\ref{prop:omegagen}). The structure of torsion sheaves on the
general twistor fibre is described by the following result.

\begin{cor}\label{cor:omegagen} Let $\pi:\kx\to\Spf(R)$ be as in Proposition
\ref{prop:omegagen}.

{\rm i)} Any torsion sheaf $F\in\coh(\kx_K)$ can be written as a
direct sum $\bigoplus F_i$ such that each $F_i$ admits a
filtration with all quotients of the form $K(\yy_i)$ for some
point $\yy_i\in\kx_K$.

{\rm ii)} If $ F\in\coh(\kx_K)$ is a non-trivial torsion free object and
$0\ne F'\in\coh(\kx_K)$ is torsion, then $\Hom_{\kx_K}(F,F')\ne0$.
\end{cor}

\begin{proof} i) Indeed, if one lifts $F$ to an $R$-flat sheaf $E$, then
$E$ is supported on a finite union of irreducible multisections
$\kz_i\subset \kx$. If only one $\kz_1$ occurs, $E$ can be
filtered such that all quotients are isomorphic to $\ko_{\kz_1}$
which induces the claimed filtration of $F=E_K$. Thus, it suffices
to show that for two distinct multisections
$\kz_1,\kz_2\subset\kx$ inducing points $\yy_1\ne\yy_2\in\kx_K$ in
the general fibre there are no non-trivial extensions, i.e.\
$\Ext^1_{\kx_K}(K(\yy_1),K(\yy_2))=0$. If $\kz_1$ and $\kz_2$ specialize
to distinct points $y_1\ne y_2\in X$ (with multiplicities), then
this obvious by semi-continuity \eqref{eq:semicont}. For $y_1=y_2$ one still has
$\chi_K(K(\yy_1),K(\yy_2))=\chi_0(k(y_1),k(y_2))=0$ due to
\eqref{eq:chi}. Using Serre duality, it therefore
suffices to show that $\Hom_{\kx_K}(K(\yy_1),K(\yy_2))=0$. If
$f:\ko_{\kz_1}\to\ko_{\kz_2}$ is  a non-trivial homomorphism, then
its image would be the structure sheaf of a subscheme of $\kx$
contained in $\kz_1$ and in $\kz_2$. Clearly, the  irreducible
multisections $\kz_i$ do not contain any proper subschemes.

ii) By i) it suffices to show that $\Hom_{\kx_K}(F,K(\yy))\ne0$ for any
closed point $\yy\in\kx_K$ and any torsion free $F\in\coh(\kx_K)$.
Using Serre duality, one knows $\Ext^2_{\kx_K}(F,K(\yy))\cong\Hom_{\kx_K}(K(\yy),F)^*=0$. Thus,
$\chi_K(F,K(\yy))=r\cdot d>0$, where $r$ is given by
$v(\res(F))=(r,s)$ and $d$ is the degree of the multisection
corresponding to $\yy\in\kx_K$, implies the assertion.
\end{proof}

Clearly, in the decomposition $F\cong\bigoplus F_i$ we may assume
that the points $\yy_i$ are pairwise distinct, which we will
usually do.

\begin{remark}\label{rem:weaker}
Later we shall use Proposition \ref{prop:omegagen} and Corollary
\ref{cor:omegagen} under slightly weaker assumptions. One easily
checks that  it suffices to assume that the first order
neighbourhood of $\kx\to\Spf(R)$ is induced by a generic twistor
space. In fact, the only assumption that is really needed is that
$\ko_\kx$ is the only line bundle on $\kx$.
\end{remark}


\subsection{Spherical objects on the very general twistor fibre}\label{sect:nosph}

The proof of the following proposition is almost a word by word
copy of the proof of \cite[Lemma 4.1]{HMS} and is included only to
show that indeed the techniques well-known for classical K3
surfaces work as well for the general twistor fibre.

\begin{prop}\label{prop:sphgen} Suppose $\pi:\kx\to \Spf(R)$ is the formal twistor space of
a K3 surface $X$ associated to a very
general K\"ahler class. Then

{\rm i)} The structure sheaf $\ko_{\kx_K}$ is the only
indecomposable rigid object in $\coh(\kx_K)$.

{\rm ii)} Up to shift, $\ko_{\kx_K}$ is the only quasi-spherical
object in $\Db(\kx_K)$ (see \cite[Def.\ 2.5]{HMS}).

{\rm iii)} The K3 category $\Db(\kx_K)$ satisfies  condition ${\rm
(\ast)}$ in \cite[Rem.\ 2.8]{HMS}.
\end{prop}

\begin{proof} Proposition 2.14 in \cite{HMS} shows that i) implies
ii) and iii). Thus, only i) needs a proof.

First, let us show that any rigid  $F\in\coh(\kx_K)$ is torsion
free. If not, the standard exact sequence $0\to F_\mathrm{tor}\to
F\to F'\to 0$ (see Section \ref{Sect:Tor}) together with
$\Hom_{\kx}(F_\mathrm{tor},F')=0$ and \cite[Lemma 2.7]{HMS} would show
that also $F_\mathrm{tor}$ is rigid. However, due to Corollary
\ref{cor:omegagen}, i) $[F_\mathrm{tor}]\in K(\kx_K)$ equals a
direct sum of sheaves of the form $K(\yy)$. As
$\chi_K(K(\yy_1),K(\yy_2))=0$ for arbitrary points
$\yy_1,\yy_2\in\kx_K$, one also has $\chi_K(F_\mathrm{tor},F_\mathrm{tor})=0$, which obviously contradicts rigidity of a non-trivial
$F_\mathrm{tor}\in\coh(\kx_K)$.

As an illustration of the techniques, let us next prove that
$\ko_{\kx_K}$ is the only spherical object in $\Db(\kx_K)$ that is
contained in  $\coh(\kx_K)$. Suppose $F\in\coh(\kx_K)$ is
spherical and let $E\in\coh(\kx)$ be an $R$-flat torsion free lift
of $F$. Then, by \eqref{eq:chi}, one has
$2=\chi_K(F,F)=\chi_0(E_0,E_0)$, i.e.\
$v(E_0)=v(\res(F))=\pm(1,0,1)$. As $F$ (and hence $E_0$) is a
sheaf, we must have $v(E_0)=(1,0,1)$. In other words, $F$ and
$\ko_{\kx_K}$ are numerically equivalent and, in particular,
$\chi_K(\ko_{\kx_K},F)=2$. The latter implies the existence of a
non-trivial $f:F\to \ko_{\kx_K}$ or a non-trivial
$f:\ko_{\kx_K}\to F$. Now we conclude by observing that any
non-trivial $f:G_1\to G_2$ in $\coh(\kx_K)$ between torsion free
$G_1$ and  $G_2$ with $v(\res(G_1))=v(\res(G_2))=(1,0,1)$ is
necessarily an isomorphism. Indeed, kernel and image of such an
$f$ are either trivial or torsion free of rank one. Since the rank
is additive and $f\ne0$, in fact $f$ is injective. The cokernel of
the injective $f:G_1\hookrightarrow G_2$ would be an
$H\in\coh(\kx_K)$ with trivial Mukai vector $v(\res(H))=0$ and
hence $H=0$ (see Example \ref{exas:MukaiVect}, iii)), i.e.\ $f$ is
an isomorphism.

Consider now an arbitrary rigid indecomposable $F\in\coh(\kx_K)$
and let $(r_0,s_0)=v(\res(F))$. Then $\chi_K(F,F)=2r_0s_0>0$ and hence
$s_0>0$. Therefore, $\chi_K(\ko_{\kx_K},F)=r_0+s_0>0$. Suppose
$\Hom_{\kx_K}(\ko_{\kx_K},F)\ne0$ and consider a short exact sequence of
the form
\begin{equation}\label{eqn:rec}\xymatrix{0\ar[r]&
\ko_{\kx_K}^{\oplus r}\ar[r]& F\ar[r]^\xi& F'\ar[r]&0.}
\end{equation} We claim that then $F'$ must be torsion free. If
not, the extension $$\xymatrix{0\ar[r]& \ko_{\kx_K}^{\oplus
r}\ar[r]&\xi^{-1}(F'_\mathrm{tor})\ar[r]&F'_\mathrm{tor}\ar[r]& 0}$$
would necessarily be non-trivial, for $\xi^{-1}(F'_\mathrm{tor})\subset F$ is torsion free. On the other hand, by Serre duality and Corollary \ref{cor:omegagen}, i), one has $\Ext^1_{\kx_K}(F'_\mathrm{tor},\ko_{\kx_K})=0$. Indeed,
$\Ext_{\kx_K}^1(K(\yy),\ko_{\kx_K})\cong\Ext_{\kx_K}^1(\ko_{\kx_K},K(\yy))^*\cong
(R^1\pi_*\ko_\kz\otimes K)^*=0$, where $\kz\subset\kx$ is the
multisection corresponding to $\yy\in\kx_K$.

Now choose $r$ maximal in \eqref{eqn:rec}. As any $0\ne
s\in\Hom_{\kx_K}(\ko_{\kx_K},F')$  defines an injection (use $F'$
torsion free), the lift of $s$ to $\tilde s\in\Hom_{\kx_K}(\ko_{\kx_K},F)$,
which exists as $\ko_{\kx_K}$ is spherical, together with the
given inclusion $\ko_{\kx_K}^{\oplus r}\subset F$ would yield an
inclusion $\ko_{\kx_K}^{\oplus r+1}\subset F$ contradicting the
maximality of $r$. Thus, for maximal $r$ the cokernel $F'$
satisfies $\Hom_{\kx_K}(\ko_{\kx_K},F')=0$ and by \cite[Lemma 2.7]{HMS}
$F'$ would be rigid as well. Then by induction on the rank, we may
assume $F'\cong\ko_{\kx_K}^{\oplus s'}$ and, since $\ko_{\kx_K}$
is spherical, this contradicts the assumption that $F$ is
indecomposable.

Eventually, one has to deal with the case that
$\Hom_{\kx_K}(\ko_{\kx_K},F)=0$, but $\Hom_{\kx_K}(F,\ko_{\kx_K})\ne0$. To
this end, consider the reflexive hull $F\ddual$. By definition
$F\ddual=(E\ddual)_K$, where $E$ is an $R$-flat lift of $F$. As we
have seen above, $F$ and hence $E$ is torsion free. Thus, $F\to
F\ddual$ is injective. The quotient map $F\ddual\to (F\ddual/F)$
deforms non trivially if $(F\ddual/F)\ne0$, e.g.\ by deforming the
support of the quotient (use Corollary \ref{cor:omegagen}, ii)).
This would contradict the rigidity of $F$. Hence $F\cong F\ddual$.
Then $\Hom_{\kx_K}(F,\ko_{\kx_K})=\Hom_{\kx_K}(\ko_{\kx_K},F\dual)$ and we can
apply the previous discussion to the rigid sheaf $F\dual$.
\end{proof}

Let us now consider the spherical twist
$$T_K:=T_{\ko_{\kx_K}}=\Phi_{(\mathcal{I}_\Delta)_K[1]}:\Db(\kx_K)\congpf\Db(\kx_K)$$ associated to the
spherical object $\ko_{\kx_K}\in\Db(\kx_K)$, i.e.\ the Fourier--Mukai equivalence with kernel $(\mathcal{I}_\Delta)_K[1]$, where $\Delta$ is the diagonal in $\kx\times_R\kx$. We have the following consequence of the previous result, which will be used in the proof of
Proposition \ref{prop:bijratpoints}.

\begin{cor}{\bf (\cite{HMS}, Proposition 2.18.)}\label{cor:recall}
Suppose $\sigma$ is a stability condition on $\Db(\kx_K)$. If
$F\in\Db(\kx_K)$ is semi-rigid with
$\sum\dim_K\Ext^i_{\kx_K}(\ko_{\kx_K},F)=1$, then there exists an
integer $n$ such that $T_K^n(F)$ is $\sigma$-stable.
\end{cor}

For the notion of stability conditions in the sense of Bridgeland
and Douglas see \cite{B}.

\subsection{Stability conditions on the very general twistor
fibre} The next task consists of actually constructing one
explicit stability condition. Following the arguments in
\cite{HMS}, it should be possible to classify all stability
conditions on $\Db(\kx_K)$ for $\kx_K$ the very general twistor
fibre as before. However,  for our purpose this is not needed.

We shall next mimic the definition of a particular stability
condition for general non-projective K3 surfaces introduced in
\cite[Sect.\ 4]{HMS}. Fix a real number $u<-1$ and let
$\kf,\kt\subset\coh(\kx_K)$ be the full additive subcategories of
all torsion free respectively torsion sheaves $F\in \coh(\kx_K)$.

\begin{lem}\label{lem:torsiontheory}
The full subcategories $\kf,\kt \subset\coh(\kx_K)$ form a torsion
theory for the abelian category $\coh(\kx_K)$.
\end{lem}

\begin{proof}
For the definition of torsion theories see e.g.\ \cite{B}. Let
$F\in\coh(\kx_K)$ and $E\in\coh(\kx)$ with $E_K\cong F$. Consider
the short exact sequence $0\to E_\mathrm{tor}\to E\to
E/E_\mathrm{tor}\to 0$ of coherent sheaves on $\kx$. Its
restriction to $\kx_K$, i.e.\ its image in $\coh(\kx_K)$, is still
a short exact sequence, which decomposes $F$ into the torsion part
$(E_\mathrm{tor})_K$ and its torsion free part
$(E/E_\mathrm{tor})_K$. As there are no non-trivial homomorphisms
from a torsion sheaf on $\kx$ to a torsion free one, the same
holds true in $\coh(\kx_K)$.
\end{proof}

The heart of the $t$-structure associated to this torsion theory
is the  abelian category
$$\ka\subset\Db(\coh(\kx_K))\cong\Db(\kx_K)$$
consisting of all complexes $F\in\Db(\coh(\kx_K))$ concentrated in
degree $0$ and $-1$ with $\kh^0(F)\in\coh(\kx_K)$ torsion and
$\kh^{-1}(F)\in\coh(\kx_K)$ torsion free. On this heart, one
defines the additive function
$$\xymatrix{Z:\ka\ar[r]&\CC},~\xymatrix{F\ar@{|->}[r]&-u\cdot
r-s},$$ where $(r,s)=v(\res(F))$. Note that by definition $Z$
takes values only in $\RR$.

\begin{prop}\label{prop:semirig}
The above construction defines a locally finite stability
condition $\sigma$ on $\Db(\kx_K)$. Moreover, if $F\in\Db(\kx_K)$
is $\sigma$-stable and semi-rigid with $\End_{\kx_K}(F)\cong K$, then, up
to shift, $F$ is a $K$-rational point $K(\xx)$.
\end{prop}

\begin{proof}
Let us first show that $Z(F)\in\RR_{<0}$ for any non-trivial
$F\in\ka$.  The Mukai vector of a torsion $F\in\coh(\kx_K)$ is of
the form $v(\res(F))=(0,s)$ with $s=-{\rm c_2}(E_0)$, where $E$ is
an $R$-flat lift of $F$. Thus, $E_0$ is a non-trivial torsion
sheaf on $X$ with zero-dimensional support and therefore
$0<h^0(E_0)=\chi(E_0)=-{\rm c}_2(E_0)$. Hence, $Z(F)\in\RR_{<0}$.

Let now $F\in\coh(\kx_K)$ be torsion free. Then due to Lemma
\ref{lem-Langton}, there exists an $R$-flat lift $E\in\coh(\kx)$
with $E_0$ a torsion free sheaf. We have to show that $u\cdot
\rk(E_0)+s(E_0)<0$ or, equivalently, ${\rm c}_2(E_0)>\rk(E_0)\cdot
(u+1)$. The inequality is linear in short exact sequences and
holds for all ideal sheaves ${\mathcal I}_Z\subset\ko_X$ of
(possibly empty) zero-dimensional subschemes $Z\subset X$. By
induction on the rank, we can therefore reduce to the case that
$H^0(X,E_0)=0$ and $H^2(X, E_0)\cong\Hom_X(E_0,\ko_X)^*=0$. (Indeed, e.g.\ any non trivial global
section of $E_0$ induces an injection $\ko_X\subset E_0$ whose
cokernel has fewer global sections. Since $\ko_X$ clearly
satisfies the inequality, which is linear in short exact
sequences, it then suffices to verify the inequality for the
quotient which is of smaller rank than the original $E_0$. Since
torsion sheaves also satisfy the inequality, one can divide out by
the torsion part of the quotient to get again a torsion free
sheaf.) But then $\chi(E_0)\leq0$ and hence the Riemann--Roch formula yields ${\rm c}_2(E_0)\geq 2 \rk(E_0)>0$.

In order to verify the Harder--Narasimhan property of $\sigma$,
one shows that the abelian category $\ka$ is noetherian and artinian. At
the same time, this then proves that $\sigma$ is locally finite.
If $F\in\ka$, then $\kh^{-1}(F)[1],\kh^0(F)\in\ka$ and the
distinguished triangle $0\to\kh^{-1}(F)[1]\to F\to\kh^0(F)\to0$ is
thus an exact sequence in $\ka$. So, if $F\supset F_1\supset
F_2\supset\ldots$ is a descending sequence in $\ka$, then the
$\kh^{-1}$ of it form a descending sequence of torsion free
sheaves. Due to rank considerations this eventually stabilizes
(the quotients $\kh^{-1}(F_i)/\kh^{-1}(F_{i+1})$ are also torsion
free!) and from then on one has a decreasing sequence of torsion
sheaves $\kh^0(F_i)\supset\kh^0(F_{i+1})\supset\ldots$. After
choosing $R$-flat lifts and restricting to the special fibre, this
yields a decreasing filtration of sheaves on $X$ concentrated in
dimension zero, which stabilizes as well. Thus, $\ka$ is artinian.
The proof that $\ka$ is noetherian is similar.

Similar arguments also prove that $\ko_{\kx_K}[1]$ is a minimal
object (i.e.\
an object without proper subobjects) in $\ka$ and therefore $\sigma$-stable of phase one.
Indeed, if $0\to F\to\ko_{\kx_K}[1]\to G\to0$ is a decomposition
in $\ka$, then the long cohomology sequence shows $\kh^0(G)=0$ and
$\rk(\kh^{-1}(F))+\rk(\kh^{-1}(G))=1$. Hence $\rk(\kh^{-1}(G))=0$,
which would yield $G=0$, or $\rk(\kh^{-1}(F))=0$. The latter would
result in a short exact sequence $0\to\ko_{\kx_K}\to\kh^{-1}(G)\to
\kh^{0}(F)\to0$ in $\coh(\kx_K)$ with $\kh^{-1}(G)$ torsion free
of rank one and $\kh^0(F)$ torsion. As shown before, a torsion
sheaf $\kh^0(F)$ has Mukai vector $(0,s)$ with $s\geq0$ and the
Mukai vector of the torsion free rank one sheaf $\kh^{-1}(G)$ is
of the form $(1,s')$ with $s'\leq1$. The additivity of the Mukai
vector leaves only the possibility $s=0$, which implies $F=0$ (see
Example \ref{exas:MukaiVect}, iii)). Contradiction.

Suppose $F\in\ka$ is a semi-rigid stable object with
$\End_{\kx_K}(F)\cong K$. If we choose a lift $E\in\Db(\kx)$ of $F$ and
denote the Mukai vector of $E_0$ by $(r,s)$, then
$0=\chi_K(F,F)=\chi(E_0,E_0)=2rs$ \eqref{eq:chi}. Hence $r=0$ or $s=0$. On the other hand,
$\chi_K(\ko_{\kx_K},F)=r+s$ and since $\ko_{\kx_K}[1],F\in\ka$ are
both non-isomorphic stable objects of the same phase,
$\Ext_{\kx_K}^i(\ko_{\kx_K},F)\ne 0$ at most for $i=0$ (use Serre
duality for $i>1$). This shows $r+s\geq0$. Thus, if $s=0$, then
$r\geq0$ and hence $r=0$, because objects in $\ka$ have
non-positive rank. Therefore, any semi-rigid stable $F\in\ka$ with
$\End_{\kx_K}(F)\cong K$ satisfies $r=0$, i.e.\ $F\in\coh(\kx_K)$, and,
moreover, $F$ is torsion. Pick an $R$-flat lift $E$ of $F$, which
is necessarily torsion as well. Proposition \ref{prop:omegagen}
shows that the support $\kz\subset\kx$ of $E$ is of relative
dimension zero over $\Spf(R)$. Clearly, the support of $E$ is
irreducible, as otherwise $F$ would have a proper subsheaf
contradicting the stability of $F$. The same argument shows that
$E$ is a rank one sheaf on $\kz$. Hence, $F\cong(\ko_\kz)_K$,
which is $K$-rational if and only if $\kz\subset\kx$ is a section
of $\pi:\kx\to\Spf(R)$. Hence $F\cong K(\xx)$ with
$\End_{\kx_K}(K(\xx))=K$.
\end{proof}


\subsection{Derived equivalences of the very general twistor fibre}\label{sect:verygentwistor}

Let us now consider two K3
surfaces $X$ and $X'$, and formal deformations of them
\begin{equation*}\xymatrix{\pi:\kx\ar[r]&\Spf(R){\text~~~~~~~~~
}{\rm and}{\text~~~~~~~~~}\pi':\kx'\ar[r]&\Spf(R).}
\end{equation*} Moreover, we shall assume that $\pi:\kx\to\Spf(R)$
is the formal twistor space of $X$ associated to a very general K\"ahler
class $\omega$.

The aim of this section is to show that under the genericity
assumption on the K\"ahler class any Fourier--Mukai equivalence
between the general fibres $\kx_K$ and $\kx_K'$ of the two formal
deformations has, up to shift and spherical twist, a sheaf kernel.

\begin{prop}\label{prop:bijratpoints}
Suppose $\Phi:\Db(\kx_K)\congpf\Db(\kx'_K)$ is a $K$-linear exact
equivalence. Then, up to shift and spherical twist in
$\ko_{\kx'_K}$, the equivalence $\Phi$ identifies $K$-rational
points of $\kx_K$ with $K$-rational points of $\kx'_K$. More
precisely, there exist integers $n$ and $m$ such that
$$T_K^n\circ\Phi[m]:\{K(\xx)~|~\xx\in\kx_K;~\End_{\kx_K}(K(\xx))\cong K\}\congpf
\{K(\yy)~|~\yy\in\kx'_K;~\End_{\kx_K}(K(\yy))\cong K\}.$$
\end{prop}

\begin{proof}
First note that although $\kx_K'$ is not necessarily a very
general twistor fibre, Proposition \ref{prop:sphgen} still applies. Indeed, the object $\ko_{\kx'_K}$ is spherical and by Proposition
2.14 we know that up to shift  $\ko_{\kx_K}$ is the only
quasi-spherical object in $\Db(\kx_K)\cong\Db(\kx'_K)$. Since the
property of being (quasi-)spherical is invariant under
equivalence, one concludes that $\ko_{\kx_K}$ is mapped to a shift
of $\ko_{\kx'_K}$ and that $\ko_{\kx'_K}$ is the only
quasi-spherical object on $\Db(\kx'_K)$.

The argument follows literally the proof of \cite[Lemma 4.9]{HMS},
so we will be brief. Let $\tilde\sigma$ be the stability condition on $\Db(\kx'_K)$ which is the image of $\sigma$
under $\Phi$. Then there exists an integer $n$, such that all
sections $K(\yy)$ of $\kx'\to\Spf(R)$ are
$T_K^n(\tilde\sigma)$-stable (cf.\  Corollary \ref{cor:recall} and
\cite[Prop.\ 2.18, Cor.\ 2.19]{HMS}). In other words, for any
$K$-rational point $\yy\in\kx'_K$ the object
$\Phi^{-1}T_K^{-n}(K(\yy))$ is $\sigma$-stable. As $\Phi$ is an
equivalence, $\Phi^{-1}T_K^{-n}(K(\yy))$ is semi-rigid
as well with $\End_{\kx_K}(\Phi^{-1}T_K^{-n}(K(\yy)))=K$. Hence, by Proposition
\ref{prop:semirig} the set $\{K(\yy)\}$  is contained in
$\{T_K^n\Phi(K(\xx))[m]\}$ for some $m$.

Applying the same argument to $\Phi^{-1}$ yields equality
$\{T^n\Phi(K(\xx)[m])\}=\{K(\yy)\}$.
\end{proof}

\begin{prop}\label{prop:kernelissheaf}
Suppose $\Phi:\Db(\kx_K)\congpf\Db(\kx'_K)$ is a Fourier--Mukai
equivalence with kernel $\ke_K\in\Db((\kx\times_R\kx')_K)$. If
$\Phi$ induces a bijection of the $K$-rational points
$\{K(\xx)\}\congpf\{K(\xx')\}$, then $\ke_K$ is a sheaf, i.e.\
$\ke_K\in\coh((\kx\times_R\kx')_K)$.
\end{prop}

\begin{proof} The full  triangulated subcategory
$\kd\subset\Db(\kx_K)\cong\Db(\coh(\kx_K))$ of all complexes
$F\in\Db(\kx_K)$ for which $\Ext_{\kx_K}^i(F,K(\xx))=0$ for all $i$ and
all $K$-rational points $\xx\in\kx_K$ will play a central role in
the proof.

 i) We shall use the following general fact: Let $F
\in\coh(\kx_K)$ such that $\Hom_{\kx_K}(F,K(\xx))=0$ for any
$K$-rational point $\xx\in\kx_K$, then $\Hom_{\kx_K}(K(\xx),F)=0$ for
any $K$-rational point $\xx\in\kx_K$. Moreover, in this case
$F\in\kd$.

In order to prove this, choose an $R$-flat lift $E$ of $F$. Then
the support of $E$ is either $\kx$ or a finite union
$\bigcup\kz_i$ of irreducible multisections (see Proposition
\ref{prop:omegagen}). In the first case we would have
$\Hom_{\kx_K}(F,K(\xx))\ne0$ for any point $\xx\in\kx_K$ (see
Corollary \ref{cor:omegagen}, ii)), contradicting the assumption.
Thus, $F\cong\bigoplus F_i$ with each $(F_i)_K$ admitting a
filtration with quotients  isomorphic to $K(\yy_i)$, where the
$\yy_i$ are points of the general fibre corresponding to different
irreducible multisections (see Corollary \ref{cor:omegagen}, i)).
By our assumption, none of the points $\yy_i$ can be $K$-rational.
But then in fact $\Hom_{\kx_K}(K(\xx),K(\yy_i))=0$ for all
$K$-rational points $\xx\in\kx_K$.

Since the $R$-flat lift $E$ of $F$ is supported in a finite union of multisection, the restriction $E_0$ of $E$ to the special fibre $X$ has rank zero. Hence $0=\rk(E_0)=\chi_0(k(x),E_0)=\chi_K(K(\xx),F)=-\dim_K\Ext^1_{\kx_K}(K(\xx),F)$, where $x\in X$ is the specialization of $K(\xx)$.
This is the second assertion.

ii) Next we claim that if $F\in\kd$, then all cohomology sheaves
$\kh^q(F)\in\coh(\kx_K)$ are as well contained in $\kd$. Indeed,
using the spectral sequence
$$E_2^{p,q}=\Ext_{\kx_K}^p(\kh^{-q}(F),K(\xx))\Rightarrow\Ext_{\kx_K}^{p+q}(F,K(\xx))$$
one sees that for $q$ minimal with non-vanishing $\kh^{-q}(F)\ne0$
any non-trivial element in
$E_2^{0,q}=\Hom_{\kx_K}(\kh^{-q}(F),K(\xx))$ would survive and
thus contradict $F\in\kd$. Hence, the maximal non-trivial
cohomology sheaf of $F$ does not admit non-trivial homomorphisms
to any $K$-rational point and is, therefore, due to i) contained
in $\kd$. Replacing $F$ by the cone of the natural morphism
$F\to\kh^{-q}(F)[q]$, which is again in $\kd$ and with a smaller number of non-trivial cohomology sheaves, one can continue and eventually proves that all
cohomology of $F$ is contained in $\kd$.

iii) Consider a  sheaf $0\ne F\in\kd\cap\coh(\kx_K)$. We claim
that $\Ext_{\kx_K}^i(\ko_{\kx_K},F)=0$ for $i\ne0$ and
$\Hom_{\kx_K}(\ko_{\kx_K},F)\ne0$.

By the definition of $\kd$, one has $\chi_K(F,K(\xx))=0$ for all
$K$-rational points $\xx\in\kx_K$. Writing this as the Mukai
pairing, one finds that the restriction $E_0$ of any $R$-flat lift
$E$ of $F$ to the special fibre $X$ will be a sheaf with Mukai
vector $(0,0,s)$, i.e.\ $E_0$ is a non-trivial sheaf concentrated
in dimension zero. If $\Ext_{\kx_K}^i(\ko_{\kx_K},F)\ne0$ for $i=1$ or
$i=2$, then, by \eqref{eq:semicont},
one would have $\Ext^i_X(\ko_{X},E_0)\ne0$, which is absurd. On
the other hand, since $s\ne0$ for $F\ne0$ (see Example
\ref{exas:MukaiVect}, iii)), the Mukai pairing also shows $0\ne
s=\chi_0(\ko_X,E_0)=\chi_K(\ko_{\kx_K},F)$ and thus
$\Hom_{\kx_K}(\ko_{\kx_K},F)\ne0$.

iv) If $F\in\kd$, then
$\Hom_{\kx_K}(\ko_{\kx_K},\kh^q(F))\cong\Ext_{\kx_K}^q(\ko_{\kx_K},F)$.
Using ii) and iii), this follows from the spectral sequence
$$E_2^{p,q}=\Ext_{\kx_K}^p(\ko_{\kx_K},\kh^q(F))\Rightarrow\Ext_{\kx_K}^{p+q}(\ko_{\kx_K},F).$$

v) Let us show that under $\Phi$ the image of any sheaf
$F\in\coh(\kx_K)$ orthogonal to all $K$-rational points is again a
sheaf $\Phi(F)\in\coh(\kx_K)$ (and moreover orthogonal to all
$K$-rational points).

As all $K$-rational points are again of the form $\Phi(K(\xx))$
for some $K$-rational point $\xx\in\kx_K$, the assumption
$F\in\kd$ implies $\Phi(F)\in\kd'$, where $\kd'\subset\Db(\kx_K')$
is defined analogously  to $\kd\subset\Db(\kx_K)$. Hence, using
iii) and iv)
$$\Hom_{\kx'_K}(\ko_{\kx'_K},\kh^q(\Phi(F)))=\Ext_{\kx'_K}^q(\ko_{\kx'_K},\Phi(F))=\Ext_{\kx_K}^q(\ko_{\kx_K},F)=0$$
for $q\ne0$. Here we use $\Phi(\ko_{\kx_K})\cong\ko_{\kx'_K}$,
which follows from Proposition \ref{prop:sphgen} saying that
$\ko_{\kx_K}\in\Db(\kx_K)$ respectively
$\ko_{\kx'_K}\in\Db(\kx'_K)$ are the only spherical objects up to
shift. (The shift is indeed trivial which follows easily from the
assumption $\Phi(K(\xx))\cong K(\xx')$.) On the other hand, by
ii), $\kh^q(\Phi(F))\in\kd'$ and thus by iii)
$\Hom_{\kx'_K}(\ko_{\kx_K'},\kh^q(\Phi(F)))\ne0$ whenever
$\kh^q(\Phi(F))\ne0$. This yields $\Phi(F)\cong\kh^0(\Phi(F))$.

vi) We will now show that $\Phi$ not only sends $K$-rational
points to $K$-rational points, but that in fact any point
$K(\yy)$, $K$-rational or not, is mapped to a point. Applying the
same argument to the inverse functor, one finds that $\Phi$
induces a bijection of the set of all ($K$-rational or not)
points.

If $K(\yy)$ is not $K$-rational, then $K(\yy)\in\kd$. Hence
$G_K:=\Phi(K(\yy))\in\kd'\cap\coh(\kx_K')$ due to v). Suppose
$G\in\coh(\kx')$ is an $R$-flat lift of $G_K$. We shall argue as
in i). Note that we can in fact apply Proposition
\ref{prop:omegagen} and Corollary \ref{cor:omegagen}, for
$\ko_{\kx'}$ is the only line bundle on $\kx'$ (otherwise there
would be an extra spherical object) and Remark \ref{rem:weaker}
therefore applies. The support of $G$ can either be $\kx'$ or a
finite union of multisections. In the first case
$\Hom_K(G_K,K(\xx'))\ne0$ for any $K$-rational point
$\xx'\in\kx'_K$. As this would contradict $G_K\in\kd'$, we
conclude that $G$ is supported on a finite union $\bigcup\kz_i$ of
multisections $\kz_i$ each inducing a point $\yy_i\in\kx_K'$.
Thus, $G_K\cong\bigoplus _{i=1}^nG_i$
 with $G_i$ admitting a filtration with quotients isomorphic to
  $K(\yy_i)$ (cf.\ Corollary \ref{cor:omegagen}, i)). Since $\Phi$ is an equivalence,
  $G_K$ is simple, i.e.\ $\End_{\kx_K}(G_K)$ is a field. Thus, $n=1$ and
$G_K=G_1\cong K(\yy_1)$.

vii) The last step is a standard argument. We have to show that
the kernel of a Fourier--Mukai transform that sends points to
points is a sheaf (cf.\ e.g.\ \cite[Lemma 3.31]{HFM}). If
$\ke\in\Db(\kx\times_R\kx')$ is a lift of $\ke_K$, we have to show
that the cohomology $\kh^q(\ke)$ for $q\ne0$ is $R$-torsion or
equivalently that $\kh^q(\ke)_K=0$ for $q\ne0$.

Suppose  $\kh^q(\ke)$ is not $R$-torsion for some $q>0$. Let $q_0$
be maximal with this property and let $\yy\in\kx_K$ be a point
corresponding to a multisection $\kz\subset\kx$ in the image of
the support of $\kh^{q_0}(\ke)$ under the first projection. Then
the sheaf pull-back $\kh^0(i^*\kh^{q_0}(\ke))$ is non-trivial,
where $i:\kz\times_R\kx'\to \kx\times_R\kx'$ is the natural
morphism. In fact, $\kh^0(i^*\kh^{q_0}(\ke))_K\ne0$. Consider the
spectral sequence (in $\Db(\kx\times_R\kx')$):
\begin{equation}\label{eqn:SS}
E_2^{p,q}=\kh^p(i^*\kh^q(\ke))\Rightarrow\kh^{p+q}(i^*\ke),
\end{equation}
which is concentrated in the region $p\leq0$.

Due to the maximality of $q_0$, the non-vanishing
$\kh^0(i^*\kh^{q_0}(\ke))_K\ne0$ implies
$\kh^{q_0}(i^*\ke)_K\ne0$. But then also
$\kh^{q_0}(\Phi(K(\yy)))=p_{\kx'_K*}(\kh^{q_0}(i^*\ke)_K)\ne0$,
which contradicts vi) if $q_0>0$.

Suppose there exists a $q_0<0$ with $\kh^{q_0}(\ke)_K\ne0$. Choose
$q_0<0$ maximal with this property and a multisection
$\yy\in\kx_K$ in the support of (the direct image under the first
projection of) $\kh^{q_0}(\ke)$. Then
$\kh^0(i^*\kh^{q_0}(\ke))_K\ne0$ and, by  using the spectral
sequence \eqref{eqn:SS} again, $\kh^{q_0}(i^*\ke)_K\ne0$. As
above, this contradicts the assumption that $\Phi(K(\yy))$ is a
sheaf.
\end{proof}


\section{Deformation of the Fourier--Mukai kernel}\label{sect:defoFM}

In this section we deal with the obstruction to
deforming the kernel of a Fourier--Mukai equivalence sideways. To this end, we need to compare
the Kodaira--Spencer classes of the two sides of the
Fourier--Mukai equivalence. Before actually showing the
triviality of the obstruction, in Section \ref{subsect:relHoch} we adapt various (known) facts about Hochschild (co)homology to our setting.


\subsection{The obstructions}\label{subs:Obstr}

Let $X$ be a smooth projective variety and let $\pi_n:\kx_n\to\Spec(R_n)$ be a scheme smooth and proper over $R_n$ such that $X\cong\kx_n\times_{R_n}\Spec(R_0)$. Assume that there exists a \emph{deformation} $\pi_{n+1}:\kx_{n+1}\to\Spec(R_{n+1})$ of $\kx_n$ to $R_{n+1}$, i.e.\ a scheme smooth and proper over $R_{n+1}$ such that $\kx_n\cong\kx_{n+1}\times_{R_{n+1}}\Spec(R_n)$. The extension class of the short exact sequence
\begin{equation*}\label{eqn:KSabso}
\xymatrix{0\ar[r]&j_{n*}\ko_X\ar[r]^{\!\!\!\!\!\!t^ndt}&\Omega_{\kx_{n+1}}|_{\kx_n}\ar[r]&\Omega_{\kx_n}\ar[r]&0}
\end{equation*}
is the \emph{(absolute) Kodaira--Spencer class} $\widetilde \kappa_n\in\Ext^1_{\kx_n}(\Omega_{\kx_n},\ko_X)$.

The \emph{(absolute) Atiyah class} of  a complex
$E_n\in\Db(\kx_n)$ is by definition the class
$\widetilde A(E_n)\in\Ext_{\kx_n}^1(E_n,E_n\otimes\Omega_{\kx_n})$
induced by the boundary map $\widetilde\alpha$ of the short exact
sequence
\begin{equation}\label{sesDiag}\xymatrix{0\ar[r]&
J_n/J_n^2\ar[r]^-{\mu_1}&\ko_{\kx_n\times
\kx_n}/J_n^2\ar[r]&\ko_{\Delta_n}\ar[r]& 0,}
\end{equation} where $J_n\subset \ko_{\kx_n\times \kx_n}$ is the ideal sheaf
of the diagonal $\widetilde\eta_n:\kx_n\congpf\Delta_n\subset
\kx_n\times \kx_n$. (Note that the fibre products are not relative
over $\Spf(R_n)$, but over $\Spec(\CC)$.) More precisely,
$$\widetilde A(E_n):E_n\to E_n\otimes\Omega_{\kx_n}[1]$$ is obtained by viewing
the boundary map  $\widetilde\alpha_n:\ko_{\Delta_n}\to
J_n/J_n^2[1]\cong\widetilde \eta_*\Omega_{\kx_n}[1]$ of
\eqref{sesDiag} as a morphism between two Fourier--Mukai kernels
and applying the induced functor transformation to $E_n$.

The Kodaira--Spencer class $\widetilde\kappa_n$ gives rise to a
morphism ${\rm
id}_{E_n}\otimes\widetilde\kappa_n:E_n\otimes\Omega_{\kx_n}[1]\to E_n\otimes
j_{n*}\ko_X[2]$, which then can be composed with the Atiyah class
$\widetilde A(E_n)$ to give a class
\begin{equation*}\label{eqn:absOBSTR}
\widetilde o(E_n):=\widetilde A(E_n)\cdot\widetilde\kappa_n\in\Ext_{\kx_n}^2(E_n,E_n\otimes (j_{n})_*\ko_X)\cong\Ext_X^2(E_0,E_0),
 \end{equation*}
where the isomorphism is given by adjunction $j_n^*\dashv (j_{n})_*$.

As it turns out, this class is the obstruction to deform $E_n$
sideways. This is

\begin{thm}\label{thm:quote}
Suppose $E_n$ is a perfect complex on $\kx_n$ and the derived
pull-back $E_0:=j_n^* E_n$ satisfies $\Ext^0_X(E_0,E_0)\cong\CC$
and $\Ext^{<0}_X(E_0,E_0)=0$. Then there exists
$E_{n+1}\in\Dp(\kx_{n+1})$ such that $i_n^*E_{n+1}\cong E_{n}$ if
and only if $\widetilde o(E_n)=0$.\qqed
\end{thm}

The proof of the theorem, which makes use of Lieblich's
obstruction in \cite{Lieblich}, can be found in \cite[App.\
C]{HMS1}. When the deformation $\kx_n$ is integrable, e.g.\
obtained as in Example \ref{ex:twistor}, the theorem  can also be
deduced from the more general results in \cite{HT}. (The
assumptions on the non-positive $\Ext$'s are not needed in
\cite{HT}.) Our main application of this theorem concerns the case
where $E_n$ is a Fourier--Mukai kernel. In this situation the
conditions on $E_0$ are easily verified.

\medskip

Let us turn to the relative versions of $\widetilde A$ and
$\widetilde\kappa_n$. The extension class
$$\kappa_{n}\in\Ext^1_{\kx_{n}}(\Omega_{\pi_{n}},\ko_{\kx_{n}})$$
of the natural short exact sequence
$$\xymatrix{0\ar[r]&\ko_{\kx_{n}}\ar[r]&\Omega_{\kx_{n+1}}|_{\kx_{n}}\ar[r]&\Omega_{\pi_{n}}\ar[r]&0,}$$
is the \emph{relative Kodaira--Spencer class} of order $n$. Here
$\Omega_{\pi_n}$ denotes the locally free sheaf of differentials
over $R_n$.

The \emph{relative Atiyah class} of a perfect complex
$E_n\in\Dp(\kx_n)$ (cf.\ Remark \ref{rmk:functors}) is the class
$$A(E_n)\in\Ext^1_{\kx_n}(E_n,E_n\otimes\Omega_{\pi_n})$$
that, as a morphism $E_n\to E_n\otimes\Omega_{\kx_n}[1]$, is
induced by the morphism of Fourier--Mukai kernels
$$\alpha_n:\ko_{\Delta_n}\to
(I_n/I_n^2)[1]\cong\eta_{n*}\Omega_{\pi_n}[1].$$ Here
$\eta_n:\kx_n\congpf\Delta_n\subset \kx_n\times_{R_n}\kx_n$ is the
relative diagonal, $I_n$ its ideal sheaf, and $\alpha_n$ is the
boundary morphism of the short exact sequence
$$\xymatrix{0\ar[r]& I_n/I_n^2\ar[r]&\ko_{\kx_n\times_{R_n}\kx_n}/I_n^2\ar[r]&\ko_{\Delta_n}\ar[r]&0.}$$

The composition of $A(E_n)$ with ${\rm id}\otimes\kappa_n[1]$ yields
the \emph{relative obstruction class}
$$o(E_n):=A(E_n)\cdot\kappa_n\in\Ext^2_{\kx_n}(E_n,E_n).$$
Moreover, the relative Atiyah class can be used to define the \emph{relative
Chern character} of a perfect complex in the usual way. (Clearly, for vector bundles or in fact arbitrary coherent sheaves
on $X_0$ this coincides by Chern--Weil theory with the usual Chern
character.) Composition in $\Dp(\kx_n)$ and exterior product
$\Omega_{\pi_n}^{\otimes i}\to\Omega_{\pi_n}^i$ allow one to form
$$\exp(A(E_n))\in\bigoplus\Ext^i_{\kx_n}(E_n,E_n\otimes\Omega^i_{\pi_n})$$
and taking the trace yields (use $E_n$ perfect)
$${\rm ch}(E_n):={\rm tr}(\exp(A(E_n)))\in\bigoplus H^i(\kx_n,\Omega^i_{\pi_n}).$$

We compare now the absolute obstruction class $\widetilde o(E_n)$
with the relative one $o(E_n)$, in the situations of Theorem
\ref{thm:quote}.

\begin{lem}\label{cor:vanishingrelrigid} Let $E_n$ be a perfect complex on $\kx_n$ as in
Theorem \ref{thm:quote} and assume furthermore that $E_0$ is
rigid, i.e.\ $\Ext^1_X(E_0,E_0)=0$. Then $E_n\in\Dp(\kx_n)$
deforms to $E_{n+1}\in\Dp(\kx_{n+1})$ if and only if the relative
obstruction class $o(E_n)\in\Ext^2_{\kx_n}(E_n,E_n)$ is trivial.
\end{lem}

\begin{proof}
The image of the Atiyah class $\widetilde A(E_n)$ under the
natural projection
$\Omega_{\kx_n}\to\Omega_{\pi_n}$ is the relative Atiyah class $A(E_n)$.
Similarly, one can compare the two Kodaira--Spencer classes
$\widetilde\kappa_n\in\Ext_{\kx_n}^1(\Omega_{\kx_n},\ko_X)$ and
$\kappa_n\in\Ext_{\kx_n}^1(\Omega_{\pi_n},\ko_\kx)$ by the
following commutative diagram
\begin{equation*}\label{eqn:KSKS}
\xymatrix{\ko_X\ar[d]\ar[r]&\Omega_{\kx_{n+1}}|_{\kx_n}\ar@{=}[d]\ar[r]&\Omega_{\kx_n}\ar[d]\ar[r]^-{\widetilde\kappa_n}&\ko_X[1]\ar[d]\\
\ko_{\kx_n}\ar[r]&\Omega_{\kx_{n+1}}|_{\kx_n}\ar[r]&\Omega_{\pi_n}\ar[r]^-{\kappa_n}&\ko_{\kx_n}[1].}
\end{equation*}
Here,  $\ko_X\to\ko_{\kx_n}\to\Omega_{\kx_{n+1}}|_{\kx_n}$ is
given by $1\mapsto t^n\mapsto t^ndt$ and
$\Omega_{\kx_n}\to\Omega_{\pi_n}$ is the natural projection. Thus
one obtains the following commutative diagram
$$\xymatrix{E_n\ar[rr]^-{\widetilde A(E_n)}\ar[drr]_{A(E_n)}
         &&E_n\otimes\Omega_{\kx_n}[1]\ar[rr]^{{\rm id}\otimes\widetilde\kappa_n}\ar[d]
                   &&E_n\otimes\ko_X[2]\ar[d]\\
            &&E_n\otimes\Omega_{\pi_n}[1]\ar[rr]_{{\rm id}\otimes\kappa_n}
                    &&E_n\otimes\ko_{\kx_n}[2].}$$

Let us consider the natural short exact sequence
$$\xymatrix{0\ar[r]&\ko_X\ar[r]^{t^n}&\ko_{\kx_n}\ar[r]&\ko_{\kx_{n-1}}\ar[r]&0.}$$
Tensoring it with $E_n$
and applying $\Hom_{\kx_n}(E_n,-)$ to it, one obtains the long exact sequence
$$\xymatrix{\ldots\ar[r]&\Ext^1_{\kx_n}(E_n,E_{n-1})\ar[r]&\Ext^2_{\kx_n}(E_n,E_0)\ar[r]^-\varphi&\Ext^2_{\kx_n}(E_n,E_n)\ar[r]&\ldots}$$ and, by the previous discussion, $\varphi(\widetilde o(E_n))=o(E_n)$.

By \eqref{eq:semicont}, $\Ext^1_X(E_0,E_0)=0$ implies $\Ext^1_{\kx_n}(E_n,E_{n-1})\cong\Ext^1_{\kx_{n-1}}(E_{n-1},E_{n-1})=0$. The claim follows.
\end{proof}


\subsection{Hochschild (co)homology and Fourier--Mukai transforms}\label{subsect:relHoch}

Let $\pi_n:\kx_n\to\Spec(R_n)$ be a smooth proper morphism with
special fibre $X=\kx_0$ of dimension $d$. We denote by
$\eta_n:\kx_n\congpf\Delta_n\subset \kx_n\times_{R_n}\kx_n$ the
relative diagonal and define the \emph{relative Hochschild
cohomology}  as the graded $R_n$-algebra
$$\HH^*(\kx_n/R_n):=\Ext^*_{\kx_n\times_{R_n}\kx_n}(\ko_{\Delta_n},\ko_{\Delta_n})\cong
\Ext^*_{\kx_n}(\eta_n^*\ko_{\Delta_n},\ko_{\kx_n}),$$ where the
second isomorphism is obtained by adjunction. Multiplication in
$\HH^*(\kx_n/R_n)$ is given by composition in
$\Db(\kx_n\times_{R_n}\kx_n)$ and the $R_n$-algebra structure is
induced by the natural map
$R_n\to\End_{\kx_n\times_{R_n}\kx_n}(\ko_{\Delta_n})$.

Similarly, one defines the \emph{relative Hochschild homology} as
$$\HH_*(\kx_n/R_n):=\Ext^{-*}_{\kx_n\times_{R_n}\kx_n}((\eta_n)_*(\omega_{\pi_n}\!\!\!\!\!\!\dual\,\,\,)[-d],\ko_{\Delta_n}),$$
which becomes a left $\HH^*(\kx_n/R_n)$-module again via
composition in $\Db(\kx_n\times_{R_n}\kx_n)$, where $\omega_{\pi_n}$ is the relative canonical bundle (see \cite{Weibel,Cal1,BF2}).
There is a natural isomorphism
$$\HH_k(\kx_n/R_n)=\Ext^{-k}_{\kx_n\times_{R_n}\kx_n}((\eta_n)_*(\omega_{\pi_n}\!\!\!\!\!\!\dual\,\,\,)[-d],\ko_{\Delta_n})
\cong\Ext^{-k}_{\kx_n}(\ko_{\kx_n},\eta_n^*\ko_{\Delta_n})={\rm
Tor}_{k}(\ko_{\kx_n},\eta_n^*\ko_{\Delta_n}).$$
This follows from
the adjunction $(\eta_{n})_!\dashv\eta_n^*$, where
$(\eta_{n})_!:\Db(\kx_n)\to\Db(\kx_n\times_{R_n}\kx_n)$ sends $E_n$
to $(\eta_{n})_*(E_n\otimes\omega_{\pi_n}\!\!\!\!\!\!\dual\,\,\,)[-d]$.

\medskip

There is a natural isomorphism in $\Db(\kx_n)$
\begin{equation}\label{HKRI}
{\bf I}:\eta_n^*\ko_{\Delta_n}\congpf\bigoplus_i\Omega_{\pi_n}^i[i]
\end{equation}
called the (relative) \emph{Hochschild--Kostant--Rosenberg} (HKR) isomorphism. The algebraic case
was studied in \cite{Yek2}, for the absolute case see also
\cite{Markarian2,Cal2}. The very general situation of an arbitrary analytic
morphism has been discussed in \cite{BF1,BF2}. A characterization of ${\bf I}$ is given in Remark \ref{rmk:AtCh}.

The HKR-isomorphism ${\bf I}$ is used in two different ways to
define isomorphisms between Hochschild (co)homology and Dolbeault
cohomology. Firstly, compose with the inverse of ${\bf I}$ in the
first variable to obtain an isomorphism
$$\Ext^*_{\kx_n}(\eta_n^*\ko_{\Delta_n},\ko_{\kx_n})\congpf\Ext^*_{\kx_n}(\bigoplus_i\Omega_{\pi_n}^i[i],\ko_{\kx_n})
\cong\bigoplus_i H^{*-i}({\kx_n},{\bigwedge}^i\kt_{\pi_n}),$$ where $\kt_{\pi_n}$ is the relative tangent bundle. This
leads to the HKR-isomorphism for Hochschild cohomology
$$I^{\rm
H\!K\!R}:\Ext^k_{\kx_n\times_{R_n}\kx_n}(\ko_{\Delta_n},\ko_{\Delta_n})=\HH^k(\kx_n/R_n)\congpf
\HT^k(\kx_n/R_n):=\bigoplus_{p+q=k}H^p(\kx_n,{\bigwedge}^q\kt_{\pi_n}).$$

For Hochschild homology we compose with ${\bf I}$ in the second
variable to obtain:
$$\Ext_{\kx_n}^{-*}(\ko_{\kx_n},\eta_n^*\ko_{\Delta_n})
\congpf\Ext_{\kx_n}^{-*}(\ko_{\kx_n},\bigoplus_i\Omega_{\pi_n}^i[i])
\cong\bigoplus_i H^{-*+i}(\kx_n,\Omega^i_{\pi_n}).$$ This leads to
the HKR-isomorphism
$$I_{\rm
H\!K\!R}:\Ext^{-k}_{\kx_n\times_{R_n}\kx_n}((\eta_n)_*(\omega_{\pi_n}\!\!\!\!\!\!\dual\,\,\,)[-d],\ko_{\Delta_n})=\HH_k(\kx_n/R_n)\congpf \HO_k(\kx_n/R_n):=\bigoplus_{q-p=k}
H^p(\kx_n,\Omega^q_{\pi_n}).$$

\begin{remark}\label{conj:CalKont}
i) If $X$ is a
smooth projective variety or a compact K\"ahler manifold, then the
deformation invariance of the Hodge numbers together with the
HKR-isomorphism $\HH_*(\kx_n/R_n)\cong\HO_*(\kx_n/R_n)$
show that the $R_n$-module $\HH_*(\kx_n/R_n)$ is free.

ii) Contraction defines on $\HO_*(\kx_n/R_n)$ the structure of
a left $\HT^*(\kx_n/R_n)$-module. The above isomorphism of
$R_n$-modules $$\left(\HH^*(\kx_n/R_n),\HH_*(\kx_n/R_n)\right)\cong \left(\HT^*(\kx_n/R_n),\HO_*(\kx_n/R_n)\right)$$ is expected to be a multiplicative
isomorphism only after correcting $I^{\rm H\! K\! R}$ by ${\rm
td}(\kt_{\pi_n})^{-1/2}$ and $I_{\rm H\! K\! R}$ by ${\rm
td}(\kt_{\pi_n})^{1/2}$. See \cite{Kont,Cal2,DTT2} for a discussion of
the absolute case and \cite{CvdB,DTT1} for the ring structure of
Hochschild cohomology.
\end{remark}

\bigskip

Take $\pi'_n:\kx_n'\to\Spec(R_n)$ another smooth proper scheme over $R_n$ of relative dimension $d'$ and set $X':=\kx'_0$. Let
$\ke_n\in\Dp(\kx_n\times_{R_n}\kx_n')$ be a perfect complex defining a Fourier--Mukai transform
$$\xymatrix{\Phi_{\ke_n}:\Db(\kx_n)\ar[r]^-\sim&\Db(\kx_n').}$$
Following \cite{Cal1}, $\Phi_{\ke_n}$ induces a graded homomorphism
$$\xymatrix{\Phi_{\ke_n}^{\HH_*}:\HH_*(\kx_n/R_n)\ar[r]&\HH_{*}(\kx_n'/R_n).}$$
For later use, we recall the construction in the relative setting.

Let $d=\dim(X)$ and $d'=\dim(X')$. For two Fourier--Mukai kernels $\ke$ and $\kf$, we denote by $\ke*\kf$ the Fourier--Mukai kernel of the composition $\Phi_\ke\circ\Phi_\kf$ (see, e.g.\ \cite[Prop.\ 5.10]{HFM}). Define the functors
\begin{equation*}
 \xymatrix{\Psi_1:\Db(\kx_n\times_{R_n}\kx_n)\ar[r]&\Db(\kx_n\times_{R_n}\kx'_n)&\Psi_2:\Db(\kx_n\times_{R_n}\kx'_n)\ar[r]&\Db(\kx'_n\times_{R_n}\kx'_n)}
 \end{equation*}
by $\Psi_1(\kg):=\ke_n*\kg$ and $\Psi_2(\kg'):=\kg'*(\ke_n)_\mathrm{R}$, where $(\ke_{n})_{\rm R}:=\ke_n\!\!\!\dual\otimes p_1^*\omega_{\pi_n}[d]$ denotes the Fourier--Mukai kernel of the right adjoint of $\Phi_{\ke_n}$ and $p_1$ is the projection (see e.g.\ \cite[Prop.\ 5.9]{HFM}).
If, as above, $\Delta_n\subset\kx_n\times_{R_n}\kx_n$ and $\Delta'_n\subset\kx'_n\times_{R_n}\kx'_n$ are the relative diagonals, then the composition $\Psi_2\circ\Psi_1$ induces a natural homomorphism
\begin{equation*}
\Ext_{\kx_n\times_{R_n}\kx_n}^{-*}((\eta_{n})_*(\omega_{\pi_n}\!\!\!\!\!\!\dual\,\,\,)[-d],\ko_{\Delta_n})\xymatrix{
\ar[r]&}
\Ext_{\kx_n'\times_{R_n}\kx_n'}^{-*}(\Psi_2(\Psi_1((\eta_{n})_*(\omega_{\pi_n}\!\!\!\!\!\!\dual\,\,\,)))[-d],\Psi_2(\Psi_1(\ko_{\Delta_n}))).
\end{equation*}
An easy computation shows that $\Psi_2(\Psi_1(\ko_{\Delta_n}))=\ke_n*(\ke_{n})_{\rm R}$ and $\Psi_2(\Psi_1((\eta_{n})_*(\omega_{\pi_n}\!\!\!\!\!\!\dual\,\,\,)))=\ke_n*(\ke_{n})_{\rm
L}\otimes p_2^*\omega_{\pi_n'}\!\!\!\!\!\!\dual\,\,\,[d-d']$, where $(\ke_{n})_{\rm L}:=\ke_n\!\!\dual\otimes p_2^*\omega_{\pi_n'}[d']$
is the Fourier--Mukai kernel of the left adjoint to $\Phi_{\ke_n}$ and $p_2$ is the projection.

The adjunction morphisms $\ke_n*(\ke_n)_{\rm R}\to\ko_{\Delta_n'}$
and $\ko_{\Delta_n'}\to\ke_n*(\ke_n)_{\rm L}$ (see \cite{Cal1})
lead to natural morphisms
$$\xi_{\rm R}:\Psi_2(\Psi_1(\ko_{\Delta_n}))\xymatrix{\ar[r]&}\ko_{\Delta_n'}~~{\rm
and}~~\xi_{\rm
L}:(\eta_{n}')_*(\omega_{\pi_n'}\!\!\!\!\!\dual\,\,)[-d']\xymatrix{\ar[r]&}
\Psi_2(\Psi_1((\eta_{n})_*(\omega_{\pi_n}\!\!\!\!\!\dual\,)))[-d]$$
(which are isomorphisms if $\Phi_{\ke_n}$ is an equivalence).
Composition with both morphisms eventually leads to
$$\xymatrix{\Phi_{\ke_n}^{\HH_*}:\HH_*(\kx_n/R_n)\ar[r]&\Ext_{\kx_n'\times_{R_n}\kx_n'}^{-*}((\eta_{n}')_*(\omega_{\pi_n'}\!\!\!\!\!\dual\,\,)[-d'],\ko_{\Delta_n'})=\HH_{*}(\kx_n'/R_n).}$$
It is easy to show, using the construction, that the action on Hochschild homology is functorial.

\medskip

If $\Phi_{\ke_n}$ is an equivalence, then it induces an isomorphism at the level of Hochschild cohomology as well.
 Hence we have a pair of isomorphisms
$$\xymatrix{\left(\Phi_{\ke_n}^{\HH^*},\Phi_{\ke_n}^{\HH_*}\right):\left(\HH^*(\kx_n/R_n),\HH_*(\kx_n/R_n)\right)\ar[r]^-\sim & \left(\HH^*(\kx'_n/R_n),\HH_*(\kx'_n/R_n)\right),}$$
with $\Phi_{\ke_n}^{\HH^*}$ an isomorphism of $R_n$-algebras and $\Phi_{\ke_n}^{\HH_*}$ compatible with the $\HH^*$-module structure on both sides (this is a simple generalization to the relative setting of e.g.\ \cite[Rmk.\ 6.3]{HFM}).

Conjugating with the HKR-isomorphisms yields also isomorphisms
$$\Phi_{\ke_n}^{\HT^*}:=I^{\rm H\! K\! R}\circ\Phi_{\ke_n}^{\HH^*}\circ
(I^{\rm H\! K\! R})^{-1}:\HT^*(\kx_n/R_n)\congpf \HT^*(\kx'_n/R_n)$$ and
$$\Phi_{\ke_n}^{\HO_*}:=I_{\rm H\! K\! R}\circ\Phi_{\ke_n}^{\HH_*}\circ
I_{\rm H\! K\! R}^{-1}:\HO_*(\kx_n/R_n)\congpf\HO_*(\kx'_n/R_n),$$ which, however, are in general not
compatible with the natural multiplicative structure of the pair
$(\HT^*,\HO_*)$. Also note that $\Phi_{\ke_n}^{\HT^*}$ and
$\Phi_{\ke_n}^{\HO_*}$ are graded, but often not bigraded.

\bigskip

Finally we recall the definition of the Hochschild versions of the
Atiyah class and the Chern character (see \cite{BF2,Cal2} and the original \cite{Markarian2}). Consider the
adjunction morphism ${\rm adj}:\ko_{\Delta_n}\to(\eta_{n})_*\eta_n^*\ko_{\Delta_n}$ as a
morphism between Fourier--Mukai kernels. The associated morphism
between the Fourier--Mukai transforms applied to
$E_n\in\Db(\kx_n)$ yields the \emph{Hochschild Atiyah class}
$$\xymatrix{A\! H(E_n):E_n\ar[r]&E_n\otimes\eta_n^*\ko_{\Delta_n}.}$$

Next one defines the \emph{Hochschild Chern character} of a
perfect complex $E_n\in\Dp(\kx_n)$ as
$${\rm ch}^{\HH_*}(E_n):={\rm tr}(A\! H(E_n))\in \Hom(\ko_{\kx_n},\eta_n^*\ko_{\Delta_n})\cong\HH_0(\kx_n/R_n).$$

\begin{remark}\label{rmk:AtCh}
i) The HKR-isomorphism ${\bf
I}$ in \eqref{HKRI} can also be described in terms of the universal relative Atiyah
class (see \cite{BF2,Cal2}). Indeed, let
$\alpha_n:\ko_{\Delta_n}\to(\eta_{n})_*\Omega_{\pi_n}[1]$ be the
universal relative Atiyah class and denote by
$$\exp(\alpha_n):\ko_{\Delta_n}\to\bigoplus_i(\eta_{n})_*\Omega^i_{\pi_n}[i]$$
its exponential. The push-forward of ${\bf I}$ composed with the
adjunction map equals $\exp(\alpha_n)$, i.e.\
\begin{equation}\label{eqn:HHexp}
\xymatrix{\exp(\alpha_n):\ko_{\Delta_n}\ar[rr]^{\rm
adj}&&(\eta_{n})_*\eta^*_n\ko_{\Delta_n}\ar[rr]^-{(\eta_{n})_*{\bf
I}}&&(\eta_{n})_*\left(\bigoplus_i\Omega^i_{\pi_n}[i]\right).}
\end{equation} In other words, under the isomorphism
$$\Hom_{\kx_n}(\eta_n^*\ko_{\Delta_n},\bigoplus_i\Omega_{\pi_n}^i[i])\cong
\Hom_{\kx_n\times_{R_n}\kx_n}(\ko_{\Delta_n},(\eta_{n})_*\left(\bigoplus_i\Omega_{\pi_n}^i[i]\right))$$
given by adjunction the HKR-isomorphism ${\bf I}$ is mapped to
$\exp(\alpha_n)$.

Using this description of ${\bf I}$, one sees that $A\!
H(E_n)$ and $A(E_n)$ are related by
\begin{equation}\label{eqn:AHA}
({\rm id}_{E_n}\otimes{\bf I})\circ A\! H(E_n)=\exp(A(E_n)).
\end{equation}
(There is a small difference between $\exp(A(E_n))$
and the Atiyah--Chern character $A\! C(E_n)=\exp(-A(E_n))$ in
\cite{BF2}, which is due to a different sign convention in the
definition of the Atiyah class. It is of no importance for our
discussion.)

ii) Taking traces of \eqref{eqn:AHA} one obtains
\begin{equation}\label{Markequ}
I_{\rm H\! K\! R}\left({\rm ch}^{\HH_*}(E_n)\right)={\rm
ch}(E_n)
\end{equation}
for all perfect complexes $E_n\in\Dp(\kx_n)$. For the absolute case see e.g.\ \cite{Cal2}.
\end{remark}

Consider a perfect complex $E_n\in\Dp(\kx_n)$ as a
Fourier--Mukai kernel on $\kx_n\cong\Spec(R_n)\times_{R_n}\kx_n$.
The induced Fourier--Mukai transform
$\Phi_{E_n}:\Db(\Spec(R_n))\to\Db(\kx_n)$ yields
$$\xymatrix{\Phi_{E_n}^{\HH_*}:\HH_0(\Spec(R_n)/R_n)\ar[r]&\HH_0(\kx_n/R_n),}$$
where, by definition, $\HH_0(\Spec(R_n)/R_n)=\End_{R_n}(R_n)$.

\begin{lem}\label{lem:altCh}
Let $\Phi_{\ke_n}:\Db(\kx_n)\to\Db(\kx_n')$ be a Fourier--Mukai functor with $\ke_n\in\Dp(\kx_n\times_{R_n}\kx'_n)$ and let $E_n\in\Dp(\kx_n)$. Then

{\rm i)} ${\rm ch}^{\HH_*}(E_n)=\Phi_{E_n}^{\HH_*}(1)$;

{\rm ii)} $\Phi_{\ke_n}^{\HH_*}\left({\rm ch}^{\HH_*}(E_n)\right)={\rm
ch}^{\HH_*}\left(\Phi_{\ke_n}(E_n)\right)$.
\end{lem}


\begin{proof}
To prove part i), one first observes that
$\Psi_2\circ\Psi_1:\End_{R_n}(R_n)\to\End_{R_n}(\Psi_2( E_n))$ sends
$1\in\End_{R_n}(R_n)$ to ${\rm id}_{\Psi_2( E_n)}$. Thus we have to show that $\xi_{\rm R}\circ\xi_{\rm
L}:(\eta_{n})_*(\omega_{\pi_n}\!\!\!\!\!\!\dual\,\,\,)[-d]\to E_n\!\!\dual\,\boxtimes E_n\to\ko_{\Delta_n}$
equals ${\rm tr}(A\! H( E_n))$ under the adjunction
$$\Hom_{\kx_n\times_{R_n}\kx_n}((\eta_{n})_*(\omega_{\pi_n}\!\!\!\!\!\!\dual\,\,\,)[-d],\ko_{\kx_n})\cong\Hom_{\kx_n}(\ko_{\kx_n},\eta_n^*\ko_{\Delta_n}).$$
To see this, observe that $\xi_{\rm L}$ under the adjunction
$$\Hom_{\kx_n\times_{R_n}\kx_n}((\eta_{n})_*(\omega_{\pi_n}\!\!\!\!\!\!\dual\,\,\,)[-d], E_n\!\!\dual\,\,\boxtimes E_n)\cong
\Hom_{\kx_n}(\ko_{\kx_n},\eta_n^*( E_n\!\!\dual\,\,\boxtimes E_n))$$
corresponds to the identity section of $\ke nd( E_n)$ and that
$\eta_n^*\xi_{\rm R}\in\Hom_{\kx_n}(\ke nd( E_n),\eta_n^*\ko_{\Delta_n})$
is obtained by composing the pull-back of the natural adjunction morphism
\begin{equation}\label{eqn:adjF}
\eta_n^*\left( E_n\!\!\dual\boxtimes E_n\to(\eta_{n})_*\eta_n^*( E_n\!\!\dual\,\boxtimes E_n)\right)
\end{equation}
with $\eta_n^*(\eta_{n})_*({\rm tr}):\eta_n^*(\eta_{n})_*\ke
nd( E_n)\to\eta_n^*\ko_{\Delta_n}$. Both assertions follow from
the construction of the left and right adjoint of a Fourier--Mukai
functor. In particular, $\xi_{\rm R}$ is indeed the composition of the
restriction $ E_n\!\!\dual\,\boxtimes E_n\to(\eta_{n})_*\ke
nd( E_n)$ with $(\eta_{n})_*({\rm tr})$. To conclude, recall that
$A\! H( E_n)$ is obtained by applying ${\rm
adj}:\ko_{\Delta_n}\to(\eta_n)_*\eta_n^*\ko_{\Delta_n}$, viewed as a
morphism of Fourier--Mukai kernels, to $ E_n$. Hence ${\rm
ch}^{\HH_*}( E_n)={\rm tr}(A\! H( E_n))$ can be described as
the composition of the identity $\ko_{\kx_n}\to\ke nd( E_n)$ with
${\rm adj}$ applied to $\ke nd( E_n)$:
\begin{equation}\label{eqn:adjF2}
\xymatrix{\ke nd( E_n)\ar[r]&\ke nd( E_n)\otimes\eta_n^*\ko_{\Delta_n}}
\end{equation}
followed by the trace on $\ke nd( E_n)$. Then note that
\eqref{eqn:adjF2} clearly equals \eqref{eqn:adjF}, which yields
the equality in i).

Part ii) is the relative version of \cite[Thm.\ 10]{Cal1} and follows directly from i) and functoriality.
\end{proof}

Notice that the right hand side in part i) of the previous lemma is taken as the definition of the Hochschild Chern character in the absolute case in \cite{Cal1}.


\subsection{Controlling the obstructions}\label{subs:FO}

Let $\kx_n$ and $\kx_n'$ be as in the previous section. Take $\ke_n\in\Dp(\kx_n\times_{R_n}\kx'_n)$ be the kernel of a Fourier--Mukai equivalence $\Phi_{\ke_n}:\Db(\kx_n)\to\Db(\kx'_n)$.
Let $v_n\in H^1(\kx_n,\kt_{\pi_n})\subset\HT^2(\kx_n/R_n)$ and
suppose $v'_n:=\Phi_{\ke_n}^{\HT^*}(v_n)\in
H^1(\kx'_n,\kt_{\pi'_n})\subset\HT^2(\kx_n/R_n)$. The inverse
images of both classes to the product $\kx_n\times_{R_n}\kx'_n$
can be considered as classes in
$H^1(\kx_n\times_{R_n}\kx'_n,\kt_{\pi_n\times\pi'_n})\subset\HT^2(\kx_n\times_{R_n}\kx'_n/R_n)$. We write $p_1^*v_n+p_2^*v_n'$ as
$v_n\boxplus v'_n\in H^1(\kx_n\times_{R_n}\kx'_n,\kt_{\pi_n\times\pi'_n})$.

For any $w\in\HT^*(\kx_n\times_{R_n}\kx'_n/R_n)$ contraction
with the exponential of the relative Atiyah class of $\ke_n$ defines a class $\exp(A(\ke_n))\cdot
w\in\Ext^*_{\kx_n\times_{R_n}\kx'_n}(\ke_n,\ke_n)$. Applied to
$v_n\boxplus v'_n$ the component in
$\Ext^2_{\kx_n\times_{R_n}\kx'_n}(\ke_n,\ke_n)$ is simply the
contraction with the relative Atiyah class.

The following is a straightforward generalization of a result of
Toda \cite{Toda} which itself relies on C\u{a}ld\u{a}raru's paper
\cite{Cal2}. For the reader's convenience we will sketch the main
arguments of the proof.

\begin{lem}\label{prop:vanishobs}
With the above assumptions, one has
$$0=A(\ke_n)\cdot(v_n\boxplus v_n')\in\Ext^2_{\kx_n\times_{R_n}\kx'_n}(\ke_n,\ke_n).$$
\end{lem}

\begin{proof}
This results relies on the commutativity of the following diagram
(see \cite[Lemma 5.7,5.8]{Toda}) which in turn is based upon the
description of the HKR-isomorphisms in terms of the universal
Atiyah class (see \eqref{eqn:HHexp}):
$$\xymatrix@C=25pt{\Ext_{\kx_n'\times_{R_n}
\kx_n'}^*(\ko_{\Delta_n'},\ko_{\Delta_n'})\ar[d]_-{\psi_2^{-1}}&&\HT^*(\kx_n'/R_n)\ar[d]^-{p_2^*}\ar[ll]_-{{I'}^{-1}}\\
\Ext^*_{\kx_n\times_{R_n} \kx_n'}(\ke_n,\ke_n)&&\HT^*(\kx_n\times_{R_n}
\kx_n'/R_n).\ar[ll]_--{\exp(A(\ke_n))\cdot}}$$

(Here, $\psi_2$ is induced by the equivalence $\Psi_2$ in Section \ref{subsect:relHoch}. Similarly,
$\psi_1$ further below is induced by $\Psi_1$.) A similar diagram holds true with $\kx_n'$
replaced by $\kx_n$ and pullling back via the first projection.
Then one has to add the pull-back $\tau^*$ of the automorphism
$\tau$ of $\kx_n\times_{R_n} \kx_n$ interchanging the two factors.
As we will see, the appearance of $\tau^*$ is crucial. One obtains
the commutative diagram
$$\xymatrix@C=25pt{\HT^*(\kx_n/R_n)\ar[d]_-{p_1^*}\ar[rr]^-{\tau^*\circ
I^{-1}}&&\Ext_{\kx_n\times_{R_n} \kx_n}^*(\ko_{\Delta_n},\ko_{\Delta_n})\ar[d]^-{\psi_1}\\
\HT^*(\kx_n\times_{R_n}
\kx_n')\ar[rr]^-{\exp(A(\ke_n))\cdot}&&\Ext^*_{\kx_n\times_{R_n}
\kx_n'}(\ke_n,\ke_n).}$$ The proof given in \cite{Toda}
generalizes in a straightforward way.

The action of $\tau^*$ can be best understood via the isomorphism
$\HT^*(\kx_n/R_n)\cong\HH^*(\kx_n/R_n)=\Ext_{\kx_n\times_{R_n}
\kx_n}^*(\ko_{\Delta_n},\ko_{\Delta_n})$. It turns out that
$\tau^*$ respects the bigrading of $\HT^*(\kx_n/R_n)$ and acts
by $(-1)^q$ on $H^p(\kx_n,\bigwedge^q\kt_{\pi_n})$.

The degree two component of the product of a class in
$H^1(\kx_n\times_{R_n} \kx_n',\kt_{\pi_n\times\pi_n'})$ with
$\exp(A(\ke_n))$ is simply the contraction with the Atiyah class.
Thus one finds
\[
\begin{split}
A(\ke_n)\cdot (v_n\boxplus v_n')&=A(\ke_n)\cdot p_1^*v_n+A(\ke_n)\cdot p_2^*v_n'\\&=\psi_1(\tau^*(I^{-1}(v_n)))+\psi_2^{-1}({I'}^{-1}(v_n'))\\
&=-\psi_1(I^{-1}(v_n))+\psi_2^{-1}({I'}^{-1}(v_n'))\\
&=\psi_2^{-1}I'^{-1}\left(-I'(\psi_2(\psi_1(I^{-1}(v_n))))+v_n'\right)\\
&=\psi_2^{-1}{I'}^{-1}\left(-\Phi_{\ke_n}^{\HT^2 }(v_n)+v_n'\right)=0
\end{split}
\]
which is precisely what we claimed.
\end{proof}

Combining this result and Lemma \ref{cor:vanishingrelrigid}, we have the following:

\begin{prop}\label{prop:def}
With the above assumptions, suppose furthermore that
\begin{itemize}
\item[{\rm i)}] $\ke_n\in\Dp(\kx_n\times_{R_n}\kx'_n)$ is rigid;

\item[{\rm ii)}] $v_n\boxplus v'_n\in H^1(\kx_n\times_{R_n}\kx'_n,\kt_{\pi_n\times\pi'_n})$ is the relative Kodaira--Spencer class $\kappa_n$ corresponding to some deformation $\ky_{n+1}\to\Spec(R_{n+1})$ of $\kx_n\times_{R_n}\kx'_n\to\Spec(R_n)$.
\end{itemize}
Then there exists a perfect complex $\ke_{n+1}\in\Dp(\ky_{n+1})$ such that $\ke_n\cong Li_n^*\ke_{n+1}$.
\end{prop}

In the next example we will see an explicit case in which assumption ii) of the previous result can be controlled.

\begin{ex}\label{ex:KSdef}
Let $\pi:\XX\to D$ be a smooth proper analytic family over a $1$-dimensional disk $D$ with distinguished fibre $X:=\pi^{-1}(0)$. We assume that the Kodaira--Spencer map, i.e.\ the boundary map of the tangent bundle sequence $\xymatrix{0\ar[r]&\kt_\pi\ar[r]&\kt_\XX\ar[r]&\pi^*\kt_D\ar[r]&0,}$
yields an isomorphism $\kappa:\kt_D\congpf R^1\pi_*\kt_\pi\cong\ke xt^1_\pi(\Omega_\pi,\ko_\XX)$ and that $h^1(\XX_t,\kt_{\XX_t})$ is constant. These assumptions are satisfied when $X$ is either a K3 surface or a product of K3 surfaces.

Fix $n>0$ and fix an embedding $\Spec(R_n)\subset D$ choosing a local parameter $t$ around $0\in D$ and define $\pi_n:=\pi|_{\kx_n}:\kx_n:=\XX\times_D \Spec(R_n)\to\Spec(R_n)$. Let $\kappa_{n-1}\in\Ext^1_{\kx_{n-1}}(\Omega_{\pi_{n-1}},\ko_{\kx_{n-1}})$ be the relative Kodaira--Spencer class and suppose that $\beta\in\Ext^1_{\kx_{n}}(\Omega_{\pi_{n}},\ko_{\kx_{n}})$ is a lift of $\kappa_{n-1}$ under the natural restriction map,
$$\xymatrix{\Ext^1_{\kx_{n}}(\Omega_{\pi_{n}},\ko_{\kx_{n}})\ar[r]&\Ext^1_{\kx_{n-1}}(\Omega_{\pi_{n-1}},\ko_{\kx_{n-1}}),}~~\xymatrix{ \beta\ar@{|->}[r]&\kappa_{n-1}.}$$
Then $\Spec(R_n)\subset D$ can be extended to $\Spec(R_{n+1})\subset D$ such that $\beta=\kappa_n$, i.e.\ $\beta$ is the relative Kodaira--Spencer class on $\kx_n$ determined by $\kx_{n+1}:=\XX\times_D\Spec(R_{n+1})$. Indeed, $\beta$ considered as a section of $\Ext^1_{\kx_n}(\Omega_{\pi_n},\ko_{\kx_n})\cong\ke xt^1_\pi(\Omega_\pi,\ko_\XX)|_{\Spec(R_n)}\cong\kt_D|_{\Spec(R_n)}$ can be locally extended to a vector field on $D$. Integrating this vector field yields a smooth curve $S\subset D$ containing $\Spec(R_n)$. The image of the restriction $\kt_S|_{\Spec(R_n)}\to\ke xt^1_{\pi_n}(\Omega_{\pi_n},\ko_{\kx_n})$ of the Kodaira--Spencer map $\kt_S\to\ke xt^1_{\pi_S}(\Omega_{\pi_S},\ko_{\XX_S})$ is thus spanned by $\beta$. Choosing the embedding $\Spec(R_{n+1})$ (i.e.\ the local parameter) appropriately, one can assume that $\beta=\kappa_n$.

Later we will consider two situations. We shall start with a deformation over a smooth one-dimensional base and study the induced finite order and formal neighbourhoods. This information will be used to construct an {\it a priori} different formal deformation by describing recursively the relative Kodaira--Spencer classes of arbitrary order.
\end{ex}


\section{Deformation of derived equivalences of K3 surfaces}\label{sect:Proof}

Let $X$ and $X'$  be two projective K3 surfaces and let
$$\Phi_{\ke_0}:\Db(X)\congpf\Db(X')$$ be a Fourier--Mukai equivalence
with kernel $\ke_0\in\Db(X\times X')$. For most of Section
\ref{sect:Proof} we will only consider the case $X=X'$. In order
to distinguish both sides of the Fourier--Mukai equivalence
however, we will nevertheless use $X'$ for the right hand side.

In this section we complete (see end of Section
\ref{sect:gospecial}) the proof of our main result, which we
restate here in a different form.

\begin{thm}\label{thm:mainbis}
 Suppose $X=X'$. Then the induced Hodge isometry $\Phi^{H^*}_{\ke_0}:\widetilde
H(X,\ZZ)\congpf\widetilde H(X,\ZZ)$ satisfies
$$\Phi^{H^*}_{\ke_0}\ne(-{\rm id}_{H^2})\oplus{\rm
id}_{H^0\oplus H^4}.$$
\end{thm}

As all orientation preserving Hodge isometries do lift to
autoequivalences (see \cite{HOY,HS,P} or \cite[Ch.\ 10]{HFM}), this
seemingly weaker form is equivalent to the original Theorem 2.

The proof splits in several steps and we argue by contradiction.
First, we need to translate the hypothesis, which is in terms of
singular cohomology, into the language of Hochschild homology.
This will allow us to deform the given Fourier--Mukai kernel
sideways to first order (see Section \ref{subs:LAST}). Extending
the kernel to arbitrary order is more involved, it will take up
Section \ref{subsec:higherorder}. Using results of Lieblich, we
conclude in Section \ref{sect:gogeneral} that the Fourier--Mukai
kernel can be extended to a perfect complex on the formal scheme
and thus leads to a derived equivalence of the general fibres. The
kernel of any Fourier--Mukai equivalence of the general fibre
however has been shown in Section \ref{sect:verygentwistor} to be
a sheaf. In Section \ref{sect:gospecial} we explain how this leads
to a contradiction when going back to the special fibre.


\subsection{From singular cohomology to first order obstruction}\label{subs:LAST}

Suppose
\begin{equation}\label{eqn:beginpf}\Phi_{\ke_0}^{H^*}:\widetilde
H(X,\ZZ)\congpf\widetilde H(X',\ZZ)
\end{equation}
preserves the K\"ahler cone up to sign, i.e.\
$\Phi_{\ke_0}^{H^*}(\kk_X)=\pm\kk_{X'}$. In the situation of our
main theorem we will have $X=X'$ and $\Phi_{\ke_0}^{H^*}$ acts on
$H^2(X,\ZZ)$ by $-{\rm id}$ and thus indeed
$\Phi_{\ke_0}^{H^2}(\kk_X)=-\kk_X$.

Consider a real ample class $\omega$ on $X$, i.e.\ $\omega\in
\kk_X\cap(\Pic(X)\otimes\RR)$ and let $v_0\in H^1(X,\kt_X)$ be the
Kodaira--Spencer class of the first order deformation of $X$ given
by the twistor space $\XX(\omega)\to\PP(\omega)$ associated to the
K\"ahler class $\omega$. More precisely, up to scaling, $v_0$ maps
to $\omega$ under the isomorphism $H^1(X,\kt_X)\congpf
H^1(X,\Omega_X^1)$ induced by a fixed trivializing section
$\sigma\in H^0(X,\omega_X)$.

\begin{lem}\label{lem:vprimezero}
Under the above assumptions,
$v'_0:=\Phi_{\ke_0}^{\HT^*}(v_0)\in H^1(X',\kt_{X'})\subset\HT^2(X').$
\end{lem}

\begin{proof} Since $\omega$ can be written as a real linear combination of
integral ample classes and all isomorphisms are linear, it
suffices to prove the assertion for $\omega={\rm c}_1(L)$ with $L$
an ample line bundle. Then there exists a line bundle $L'\in\Pic(X')$ such that
$\Phi_{\ke_0}^{H^*}({\rm c}_1(L))={\rm c}_1(L')$.
We claim that then also $\Phi_{\ke_0}^{\HO_*}({\rm
c}_1(L))={\rm c}_1(L')$. Indeed, by Lemma \ref{lem:altCh}, ii), one knows
that $\Phi_{\ke_0}^{\HH_*}\circ{\rm ch}^{\HH_*}={\rm ch}^{\HH_*}\circ\Phi_{\ke_0}$, which combined with
\eqref{Markequ} yields $\Phi^{\HO_*}\circ{\rm ch}={\rm
ch}\circ\Phi_{\ke_0}$. On the other hand, for K3 surfaces one has
${\rm td}(X')^{1/2}\cdot({\rm ch}\circ\Phi_{\ke_0})=
v\circ\Phi_{\ke_0}=\Phi^{H^*}_{\ke_0}\circ v$, where $v:={\rm
td}^{1/2}\cdot{\rm ch}$ is the Mukai vector on $X$ respectively
$X'$. Thus, ${\rm td}(X')^{1/2}\cdot(\Phi_{\ke_0}^{\HO_*}\circ{\rm ch})={\rm
td}(X')^{1/2}\cdot({\rm
ch}\circ\Phi_{\ke_0})=\Phi_{\ke_0}^{H^*}\circ v$.
Since multiplication with ${\rm td}^{1/2}$ does not affect the
component of degree two, this shows that $\Phi_{\ke_0}^{H^*}({\rm
c}_1(L))={\rm c}_1(L')$ implies $\Phi_{\ke_0}^{\HO_*}({\rm
c}_1(L))={\rm c}_1(L')$.

In the next step we shall use the following:

\medskip

\noindent{\it Claim.} Suppose $\alpha\in \HH_0$ with $I(\alpha)\in H^1(\Omega)\subset
\HO_0$. Let $w\in \HH^2$ such that $w\cdot\sigma=\alpha$
and $w_0:=I(w)\in\HT^2$. Then
\begin{equation}\label{eqn:applRam}
    w_0\lrcorner\sigma =I(\alpha).
\end{equation}

\smallskip

Indeed, $w_0\lrcorner \sigma=I(w)\lrcorner I(\sigma)=({\rm
td}^{-1/2}I(w))\lrcorner({\rm td}^{1/2}I(\sigma))={\rm
td}^{1/2}I(w\cdot\sigma)=I(\alpha)$. Here we used $I(\alpha)\in
H^1(\Omega^1)$ for the second and the last equality (write down
the bidegree decomposition for $I(w)$ which a priori might have
components not contained in $H^1(\kt)$), $\sigma\in
H^0(\Omega^2)$ for the second one, and \cite[Thm.\ 1.2]{MNS} for the penultimate one.

\medskip

As ${\rm c}_1(L)\in
H^1(\Omega)$, by the previous claim, there exists $w_0\in\HT^2$ such that ${\rm c}_1(L)=w_0\lrcorner\sigma$.
Hence, if $\sigma':=\Phi_{\ke_0}^{\HO_*}(\sigma)$, we get the following sequence of equalities:
\[
\begin{split}
    \Phi^{\HO_*}({\rm c}_1(L))&=\Phi^{\HO_*}(w_0\lrcorner\sigma)=I\Phi^{\HH_*}I^{-1}(w_0\lrcorner\sigma)\\
&\stackrel{\rm (i)}=I\Phi^{\HH_*}(w\cdot\sigma)\stackrel{\rm (ii)}=I(\Phi^{\HH^*}(w)\cdot\sigma')\\
&\stackrel{\rm (iii)}=(I\Phi^{\HH^*}(w))\lrcorner\sigma',
\end{split}
\]
where (i) follows from the claim for $I(w)=w_0$, (ii) is due to the multiplicativity of $(\Phi^{\HH^*},\Phi^{\HH_*})$, and (iii) is obtained by applying again the claim to the $(1,1)$-class ${\rm
c}_1(L')=\Phi_{\ke_0}^{\HO_*}({\rm
c}_1(L))$. Thus, $(I\Phi^{\HH^*}(w))\lrcorner\sigma'\in H^{1,1}(X')$, which suffices to
conclude $v_0'\in H^1(X',\kt_{X'})$.
\end{proof}

\begin{remark}\label{rmk:Ramadoss}
    The result in Lemma \ref{lem:vprimezero} can be derived from the compatibility between the action of a Fourier--Mukai transform on Hochschild homology and the one on singular cohomology proved in \cite[Thm.\ 1.2]{MS} which, in turn, relies on \cite{Markarian2,Ramadoss}. The main result in \cite{Ramadoss} can also be used to deduce directly \eqref{eqn:applRam}.
\end{remark}

Using for example \cite{Toda}, we can then construct $\kx_1'\to\Spec(R_1)$ with Kodaira--Spencer class $v_0'\in H^1(X',\kt_{X'})$. Note that by construction $\kx'_1\to\Spec(R_1)$ depends on the actual Fourier--Mukai kernel $\ke_0$.

\begin{cor}
The Fourier--Mukai kernel $\ke_0$ extends to a perfect complex
$\ke_1\in\Dp(\kx_1\times_{R_1}\kx_1')$ inducing an equivalence
$\Phi_{\ke_1}:\Db(\kx_1)\congpf\Db(\kx_1')$.
\end{cor}

\begin{proof}
As $\ke_0$ is the Fourier--Mukai kernel of an equivalence, $\ko_\Delta\boxtimes\ke_0$ defines an equivalence $\Db(X\times X)\cong\Db(X\times X')$ that sends $\ko_\Delta$ to $\ke_0$. In particular, this yields an isomorphism $\Ext^1_{X\times X}(\ko_\Delta,\ko_\Delta)\cong\Ext^1_{X\times X'}(\ke_0,\ke_0)$. But since $\Ext^1_{X\times X}(\ko_\Delta,\ko_\Delta)\cong H^0(X,\kt_X)\oplus H^1(X,\ko_X)=0$ for the K3 surface $X$, this immediately shows that $\ke_0$ is rigid.
Hence the existence of $\ke_1$ follows from Proposition \ref{prop:def}.
The assertion that $\Phi_{\ke_1}$ is again an equivalence is part iii) of Remark \ref{rmk:functors}.
\end{proof}


\subsection{Deforming to higher order}\label{subsec:higherorder}

The idea to proceed is to extend recursively $\kx_n'$ to
$\kx_{n+1}'$ such that the Fourier--Mukai kernel $\ke_n$ on
$\kx_n\times_{R_n}\kx_n'$ deforms to a perfect complex on
$\kx_{n+1}\times_{R_{n+1}}\kx_{n+1}'$.

We will choose a formal twistor space $\pi:\kx\to\Spf(R)$ of $X$ associated to
a very general real ample class $\omega$ with $n$-th order neighbourhoods
$\pi_n:\kx_n\to\Spec(R_n)$. The
relative Kodaira--Spencer classes of $\kx_n\subset\kx_{n+1}$ will
be denoted $v_n\in H^1(\kx_n,\kt_{\pi_n})$.

Suppose we have constructed $\pi_n':\kx_n'\to\Spec(R_n)$ and a
perfect complex $\ke_n\in\Dp(\kx_n\times_{R_n}\kx_n')$ such that
$\Phi:=\Phi_{\ke_n}:\Db(\kx_n)\congpf\Db(\kx_n')$ is an
equivalence. Then let
\begin{equation}\label{eqn:defvn}
    v_n':=\Phi^{\HT^2}(v_n)\in\HT^2(\kx'_n/R_n).
\end{equation}

We would like to view $v_n'$ as a relative Kodaira--Spencer class
of order $n$ on $\kx'_n$ of some extension
$\kx_n'\subset\kx_{n+1}'\to \Spec(R_{n+1})$. For this we need the
following lemma, which is the higher order version of Lemma
\ref{lem:vprimezero}. However, the reader will observe that the
arguments in the two situations are different and neither of the
proofs can be adapted to cover the other case as well.

\begin{lem}\label{lem:cruc}
The class $v_n'$ is contained in $H^1(\kx_n',\kt_{\pi'_n})\subset
\HT^2(\kx'_n/R_n)$.
\end{lem}

\begin{proof}
Let $\sigma_n\in \HH_2(\kx_n/R_n)=H^0(\kx_n,\omega_{\pi_n})=\HO_2(\kx_n/R_n)$ be a trivializing section of $\omega_{\pi_n}$
and let $\sigma_n':=\pm\Phi^{\HO_*}(\sigma_n)\in \HH_2(\kx'_n/R_n)$. Furthermore, let
$\omega_n:=v_n\lrcorner\sigma_n\in
H^1(\kx_n,\Omega_{\pi_n})\subset\HO_0(\kx_n/R_n)$ and
$\omega_n':=v_n'\lrcorner\sigma_n'\in\HO_0(\kx'_n/R_n)$.
Then also $\omega_n'=\pm\Phi^{\HO_0}(\omega_n)$, as
$(\Phi^{\HH^*},\Phi^{\HH_*})$ is compatible with the
multiplicative structure.

Clearly, $v_n'$ is contained in $H^1(\kx_n',\kt_{\pi'_n})\subset
\HT^2(\kx_n'/R_n)$ if and only if $\omega_n'$ is contained in
$H^1(\kx_n',\Omega_{\pi_n'})\subset\HO_0(\kx_n'/R_n)$.

In a first step, we shall show that $\Phi^{\HO_0}:\HO_0(\kx_n/R_n)\congpf\HO_0(\kx_n'/R_n)$ preserves
$H^{0,0}\oplus H^{2,2}$, i.e.\ that it maps $(H^{0,0}\oplus
H^{2,2})(\kx_n/R_n):=H^0(\kx_n,\ko_{\kx_n})\oplus
H^2(\kx_n,\omega_{\pi_n})$ to $(H^{0,0}\oplus
H^{2,2})(\kx_n'/R_n):=H^0(\kx_n',\ko_{\kx'_n})\oplus
H^2(\kx_n',\omega_{\pi_n'})$.

To this end, consider the Chern character ${\rm ch}(E_n)\in\HO_0(\kx_n/R_n)$ for arbitrary
$E_n\in\Dp(\kx_n)$. In particular, ${\rm ch}(\ko_{\kx_n})=1\in
H^0(\kx_n,\ko_{\kx_n})\subset\HO_0(\kx_n/R_n)$, since
$A(\ko_{\kx_n})$ is by definition trivial. Furthermore, if
$k(x_n)\in\Db(\kx_n)$ denotes the structure sheaf of a section of
$\pi_n:\kx_n\to\Spec(R_n)$, then ${\rm ch}(k(x_n))$ is contained
in $H^2(\kx_n,\omega_{\pi_n})\subset\HO_0(\kx_n/R_n)$, as
rank and determinant of $k(x_n)$ are trivial. Actually, ${\rm
ch}(k(x_n))$ trivializes the $R_n$-module
$H^2(\kx_n,\omega_{\pi_n})$. Indeed, since the Atiyah class is
compatible with pull-back, one has $j_n^*{\rm ch}(k(x_n))={\rm
ch}(k(x_0))$ and the latter is clearly non-trivial in
$H^{2,2}(X)=\CC$.

So, $(H^{0,0}\oplus H^{2,2})(\kx_n/R_n)$ is contained in the
$R_n$-submodule of $\HO_0(\kx_n/R_n)$ spanned by the Chern
character of perfect complexes. The analogous assertion holds true
for $\kx_n'$.

As we will show now, in fact equality holds. This will later be
needed only for $\kx_n'$. So we write it down in this case.
If $E_n'\in\Dp(\kx_n')$, then ${\rm
ch}_1(E_n')\in H^1(\kx_n,\Omega_{\pi_n'})$ equals ${\rm
tr}(A(E_n'))$, which by standard arguments is simply
$A(\det(E_n'))$. The determinant $\det(E_n')$ is a line bundle
on $\kx_n'$. Therefore, it suffices to prove that any line bundle
on $\kx_n'$ is trivial, but this has been discussed already in
Section \ref{sect:Twistor}. In fact, it suffices to prove this for
$n=1$ and then $\kx_1'\to\Spec(R_1)$ is the first infinitesimal
neighbourhood of $X'$ inside its twistor space associated to the
K\"ahler class $\omega'$ (see Remark \ref{rem:weaker}).

As explained in the proof of Lemma \ref{lem:vprimezero},
\eqref{Markequ} and part ii) of Lemma \ref{lem:altCh} imply
$\Phi^{\HO_0}\circ {\rm ch}={\rm ch}\circ\Phi$. Hence
$\Phi^{\HO_0}({\rm ch}(E_n))$ is contained in
$(H^{0,0}\oplus H^{2,2})(\kx_n'/R_n)$ for any
$E_n\in\Dp(\kx_n)$ and therefore
$$\Phi^{\HO_0}\left((H^{0,0}\oplus
H^{2,2})(\kx_n/R_n)\right)\subset (H^{0,0}\oplus
H^{2,2})(\kx_n'/R_n).$$

Let now $w_n:=I^{-1}(v_n)\in \HH^2(\kx_n/R_n)$ and
$w_n':=I^{-1}(v_n')\in \HH^2(\kx_n'/R_n)$. Then by definition of
$\Phi^{\HT^*}$ we have $\Phi^{\HH^2}(w_n)=w_n'$. The
multiplicativity of $(\Phi^{\HH^*},\Phi^{\HH_*})$ and
Lemma \ref{lem:altCh}, ii) imply $\Phi^{\HH_*}(w_n\cdot{\rm ch}^{\HH_*}(E_n))=w_n'\cdot{\rm ch}^{\HH_*}(\Phi(E_n))$. So,
$w_n\cdot{\rm ch}^{\HH_*}(E_n)=0$ for all
$E_n\in\Dp(\kx_n)$ if and only if $w_n'\cdot{\rm ch}^{\HH_*}(E_n')=0$ for all $E_n'\in\Dp(\kx_n')$.

Suppose we know already that in general
\begin{equation}\label{eqn:AH}
w_n\cdot{\rm ch}^{\HH_*}(E_n)=0~~{\rm if~and~ only~ if
~}v_n\lrcorner{\rm ch}(E_n)=0
\end{equation}
and the analogous statement on $\kx_n'$. (For this assertion our
assumptions that $\kx_n\to\Spec(R_n)$ is of dimension two and that
$v_n=I(w_n)\in\HT^2(\kx_n/R_n)$, $w_n\in \HH^2(\kx_n/R_n)$
are important. See below.) Then one concludes as follows. Since
obviously $v_n\lrcorner(H^{0,0}\oplus H^{2,2})(\kx_n/R_n)=0$ and
as shown above ${\rm Im}({\rm ch})\subset(H^{0,0}\oplus
H^{2,2})(\kx_n/R_n)$, the `if' direction in \eqref{eqn:AH} would
yield $w_n\cdot{\rm ch}^{\HH_*}(E_n)=0$ for all $E_n\in\Dp(\kx_n)$ and hence $w_n'\cdot{\rm ch}^{\HH_*}(E_n')=0$ for all $E_n'\in\Dp(\kx_n')$.
As $(H^{0,0}\oplus H^{2,2})(\kx_n'/R_n)$ is actually spanned by
${\rm Im}({\rm ch})$, the `only if' direction in \eqref{eqn:AH}
then shows that $v_n'\lrcorner(H^{0,0}\oplus H^{2,2})(\kx_n'/R_n)=0$. The latter clearly means $v_n'\in H^1(\kx_n',\kt_{\pi_n'})$.

The assertion \eqref{eqn:AH} follows almost directly from
\eqref{eqn:AHA}. More precisely, in our situation, \cite[Cor.\ 5.2.3]{BF2} says that for any $v_n\in\HT^2(\kx_n/R_n)$ and any $E_n\in\Dp(\kx_n)$ the part of $v_n\lrcorner(\exp(A(\kf_n)))$ contained in $\Ext^2_{\kx_n}(E_n,E_n)$
coincides with the projection of $I^{-1}(v_n)\cdot A\! H(E_n)$
under
$\Ext^2_{\kx_n}(E_n,E_n\otimes\eta_n^*\ko_{\Delta_n})\congpf\Ext^2_{\kx_n}(E_n,E_n)$.
Taking traces on both sides yields $v_n\lrcorner{\rm
ch}(E_n)=I^{-1}(v_n)\cdot{\rm ch}^{\HH_*}(E_n)$, which then
proves \eqref{eqn:AH}.
\end{proof}

\begin{cor}\label{cor:recurs}
If $\ke_n\in\Dp(\kx_n\times_{R_n}\kx_n')$ induces an equivalences
$\Phi_{\ke_n}:\Db(\kx_n)\congpf\Db(\kx_n')$, then there exists a
deformation $\kx_{n+1}'\to\Spec(R_{n+1})$ of $\kx_n'\to\Spec(R_n)$
and a complex $\ke_{n+1}\in\Dp(\kx_{n+1}\times_{R_{n+1}}\kx_{n+1}')$ such that $\ke_n\cong Li_n^*\ke_{n+1}$. Moreover, $\ke_{n+1}$ induces an equivalence
$\Phi_{\ke_{n+1}}:\Db(\kx_{n+1})\congpf\Db(\kx_{n+1}')$.
\end{cor}

\begin{proof}
By Example \ref{ex:KSdef}, choose the extension $\kx_{n+1}'\to\Spec(R_{n+1})$ such that
its Kodaira--Spencer class $\kappa_n\in H^1(\kx_n',\kt_{\pi_n'})$
is $v_n'$ in \eqref{eqn:defvn}, which by Lemma \ref{lem:cruc} is
indeed an element in $H^1(\kx_n',\kt_{\pi_n'})$.
Since $\ke_0=Lj_n^*\ke_n$ is rigid, Proposition \ref{prop:def} allows one
to conclude the existence of a complex
$\ke_{n+1}\in\Dp(\kx_{n+1}\times_{R_{n+1}}\kx_{n+1}')$ with
$Li_n^*\ke_{n+1}\cong\ke_n$. The last assertion follows from
part iii) of Remark \ref{rmk:functors}.
\end{proof}


\subsection{Deformation to the general fibre}\label{sect:gogeneral}

Applying Corollary \ref{cor:recurs} recursively, we obtain a
formal scheme $\pi':\kx'\to\Spf(R)$ and perfect complexes
$\ke_n\in\Dp(\kx_n\times_{R_n}\kx_n')$, $n\in\NN$, inducing Fourier--Mukai
equivalences $\Phi_{\ke_n}:\Db(\kx_n)\congpf\Db(\kx_n')$ with
$Li_n^*\ke_{n+1}\cong\ke_n$ and with $\ke_0$ as given in
\eqref{eqn:beginpf}.

Now we use Lieblich's \cite[Sect.\ 3.6]{Lieblich} to conclude that
the existence of all higher order deformations is enough to
show the existence of a formal deformation of the complex.
So, there exists a complex $\ke\in\Db(\kx\times_R\kx')$ with $L\iota_n^*\ke\cong\ke_n$, for all $n\in\NN$.

\begin{remark}
Lieblich's result is far from being trivial and the proof is quite
ingenious. Of course, if a coherent sheaf lifts to any order, it
deforms by definition to a sheaf on the formal neighbourhood. For
complexes as objects in the derived categories this is a different
matter. Note that a priori one really only gets an object in
$\Db(\kx\times_R\kx')$, which is, by definition, $\Db_{\rm
coh}(\Mod{\kx\times_R\kx'})$, and not in
$\Db(\coh(\kx\times_R\kx'))$ as one could wish for.

For the convenience of the reader, let us recall Lieblich's
strategy. Instead of considering deformations of $\ke_0$ as an
object in the derived category, Lieblich shows in \cite[Prop.\
3.3.4]{Lieblich} that by replacing $\ke_0$ with a complex of
quasi-coherent injective sheaves one can work with actual
deformations of complexes, i.e.\ the differentials and objects are
deformed (flat over the base) and the restrictions to lower order
yield isomorphisms of complexes. By taking limits, one obtains a
bounded complex of \emph{ind-quasi-coherent sheaves} on the formal
scheme. Eventually, one has to show that the complex obtained in
this way, which is an object in $\Db(\Mod{\kx\times_R\kx'})$, has
coherent cohomology. This is a local statement and is addressed in
\cite[Lemma 3.6.11]{Lieblich}. Note that the main result
\cite[Prop.\ 3.6.1]{Lieblich} treats the case that the formal
scheme is given as a formal neighbourhood of an actual scheme over
$\Spec(R)$ and asserts then the existence of a perfect complex on
the scheme. In our case, the  actual scheme does not exist but
only the formal one. However, Lieblich's arguments proving the
existence of the perfect complex on the formal scheme, which is
the first step in his approach, do not use the existence of the
scheme itself.
\end{remark}

Now, by Remark \ref{rmk:functors}, iv), $\Phi_{\ke_K}:\Db(\kx_K)\congpf\Db(\kx'_K)$
is an equivalence. The Fourier--Mukai equivalence
$T_{\ko_{\kx'}}:\Db(\kx')\congpf\Db(\kx')$ with kernel
$\ki_{\Delta_{\kx'}}[1]$ `restricted'
to the special fibre is the spherical twist
$$\xymatrix{T_0:\Db(X')\ar[r]^-\sim&\Db(X')}$$ and
`restricted' to the general fibre it yields the spherical twist
$$\xymatrix{T_K:\Db(\kx_K')\ar[r]^-\sim&\Db(\kx_K').}$$ Then
Proposition \ref{prop:bijratpoints} asserts that there exist
integers $n$ and $m$ such that the composition
$T_K^n\circ\Phi_{\ke_K}[m]$ defines a bijection between the set of
$K$-rational points of $\kx_K$ and $\kx_K'$. By the discussion in
Section \ref{sect:verygentwistor} this is enough to conclude that
$T_K^n\circ\Phi_{\ke_K}[m]$ can be written as a Fourier--Mukai
transform whose kernel is a sheaf on $(\kx\times_R\kx')_K$. Note that
$n$ and $m$ must both be even. Indeed if a $K$-rational point is
sent to a $K$-rational point via $T_K^n\circ\Phi_{\ke_K}[m]$, then
its restriction to the special fibre $T_0^n\circ\Phi_{\ke_0}[m]$
preserves the Mukai vector of a point $(0,0,1)$. Now use that
$T_0^{H^*}$ sends $(0,0,1)$ to $(-1,0,0)$, that the simple shifts
acts by $-{\rm id}$, and that $\Phi_{\ke_0}^{H^*}$ preserves
$(0,0,1)$ by assumption.

The conclusion of the discussion so far is that, up to applying
shift and spherical twist, the Fourier--Mukai kernel $\ke_0$ of an
equivalence $\Phi$ with $\Phi^{H^*}=-{\rm id}_{H^2}\oplus{\rm
id}_{H^0\oplus H^4}$ deforms to a sheaf on the general fibre
$(\kx\times_R\kx')_K$, where $\kx\to\Spf(R)$ is the formal
neighbourhood of $X$ inside a very general twistor space and
$\kx'\to\Spf(R)$ was constructed recursively. One now has to show
that this leads to a contradiction.


\subsection{Return to the special fibre}\label{sect:gospecial}

Let $X$ be a smooth projective K3 surface and $\kg$ a coherent sheaf on
$X\times X$. Consider the Fourier--Mukai transform
$\Phi_\kg:\Db(X)\to\Db(X)$ with kernel $\kg$. As we make no
further assumptions on $\kg$, $\Phi_\kg$ is not necessarily an equivalence. We
shall be interested in the induced map on cohomology
$\Phi^{H^*}_\kg:\widetilde H(X,\QQ)\to\widetilde H(X,\QQ)$.

\begin{lem}\label{lem:weirdaction} For any sheaf $\kg$ on $X\times
X$ one has $\Phi_\kg^{H^*}\ne(-{\rm id}_{H^2})\oplus{\rm
id}_{H^0\oplus H^4}$.
\end{lem}

\begin{proof}
Suppose $\Phi_\kg^{H^*}=(-{\rm id}_{H^2})\oplus{\rm id}_{H^0\oplus
H^4}$. Choose an ample line bundle $L$ on $X$. Then for $n,m\gg0$
the sheaf $\kg_{n,m}:=\kg\otimes (q^*L^n\otimes p^* L^m)$ is
globally generated and $\Phi_\kg(L^n)=p_*(\kg\otimes q^* L^n)$ is
a sheaf. So, there exists a short exact sequence
$0\to\kk\to\ko_{X\times X}^N\to\kg_{n,m}\to0$. Twisting further
with $q^*L^{n'}$, $n'\gg0$, kills the higher direct images of
$\kk$ under the projection $p$, i.e.\ $R^ip_*(\kk\otimes
q^*L^{n'})=0$ for $i>0$. Thus, there exists a surjection
$\ko_X^{N'}\cong p_*(\ko_{X\times X}^N\otimes
q^*L^{n'})\twoheadrightarrow
p_*(\kg_{n+n',m})=\Phi_\kg(L^{n+n'})\otimes L^m$.

On the other hand, by assumption $v(\Phi_\kg(L^{n+n'})\otimes
L^m)=1+(m-(n+n')){\rm c}_1(L)+s$ for some $s\in H^4(X,\QQ)$. Thus,
$\Phi_\kg(L^{n+n'})\otimes L^m$ is a globally generated coherent
sheaf of rank one with first Chern class $(m-(n+n')){\rm c}_1(L)$.
It is not difficult to see that this is impossible as soon as
$m-(n+n')<0$.
\end{proof}

We leave it to the reader to formulate a similar statement for
sheaves on the product $X\times X'$ of two not necessarily
isomorphic K3 surfaces.

Consider two formal deformations $\kx\to\Spf(R)$ and
$\kx'\to\Spf(R)$ of the same algebraic K3 surface
$X=\kx_0=\kx_0'$.

\begin{cor}\label{cor:sheafcannot}
Let $\ke\in\Db(\kx\times_R\kx')$ be an object whose restriction to
the general fibre is a sheaf, i.e.\ $\ke_K\in\coh(\kx_K)$. If
$\ke_0\in\Db(X\times X)$ denotes the restriction to the special
fibre, then the  Fourier--Mukai transform
$\Phi_{\ke_0}:\Db(X)\to\Db(X)$ induces a map
$\Phi_{\ke_0}^{H^*}:\widetilde H(X,\QQ)\to \widetilde H(X,\QQ)$
different from $(-{\rm id}_{H^2})\oplus{\rm id}_{H^0\oplus H^4}$.
\end{cor}

\begin{proof}
If $\ke_K$ is a sheaf, then there exists an $R$-flat lift
$\tilde\ke\in\coh(\kx\times_R\kx')$ of $\ke_K$. Thus the complex
$\ke$ and the sheaf $\tilde\ke$ coincide on the general fibre or,
in other words, they differ by $R$-torsion complexes. In
particular, the  restrictions to the special fibres $\ke_0$ and
$\tilde\ke_0$ define the same elements in the K-group (see Remark
\ref{rem:resonK}) and therefore the same correspondence
$\Phi_{\ke_0}^{H^*}=\Phi_{\tilde\ke_0}^{H^*}:\widetilde
H(X,\ZZ)\congpf\widetilde H(X,\ZZ)$.

So, $\tilde\ke_0$ is a sheaf(!) on $X\times X$ inducing the same
map  on cohomology as the complex $\ke_0$. Now Lemma
\ref{lem:weirdaction} applies and yields the contradiction.
\end{proof}

Clearly, the corollary applies directly to our problem with $\ke$
as in Section \ref{sect:gogeneral} and which therefore contradicts
the assumption that $\Phi_{\ke_0}^{H^*}$ acts as $(-{\rm
id}_{H^2})\oplus{\rm id}_{H^0\oplus H^4}$. This concludes the
proof of Theorem \ref{thm:mainbis}.


\subsection{Derived equivalence between non-isomorphic K3 surfaces}\label{sect:noniso}

The main theorem implies that every
derived equivalence between projective K3 surfaces is orientation
preserving.

Let $X$ be an arbitrary K3 surface. Then its cohomology
$\widetilde H(X,\ZZ)$ admits a natural orientation (of the
positive directions). Indeed, if $\omega\in H^{1,1}(X)$ is any
K\"ahler class, then $1-\omega^2/2$ and $\omega$ span a positive
plane in $\widetilde H(X,\RR)$. Another positive plane orthogonal
to it is spanned by real and imaginary part of a generator
$\sigma\in H^{2,0}(X)$. Together they span a positive four-space
which is endowed with a natural orientation by choosing the base
${\rm Re}(\sigma),{\rm Im}(\sigma),1-\omega^2/2,\omega$. This
orientation does neither depend on the particular K\"ahler class
$\omega$ nor on the choice of the regular two-form $\sigma$.

If $X'$ is another K3 surface, one says that an isometry
$\widetilde H(X,\ZZ)\congpf \widetilde H(X',\ZZ)$ is orientation
preserving if the natural orientations of the four positive
directions in $\widetilde H(X,\RR)$ respectively in $\widetilde
H(X',\RR)$ coincide under the isometry.

\begin{cor}
Let $\Phi:\Db(X)\congpf\Db(X')$ be an exact equivalence between
two projective K3 surfaces. Then the induced Hodge isometry
$$\Phi^{H^*}:\widetilde H(X,\ZZ)\congpf\widetilde H(X',\ZZ)$$ is
orientation preserving.
\end{cor}

\begin{proof}
Suppose that $\Phi^{H^*}$ is an orientation reversing Hodge isometry.
Composing with $-{\rm id}_{H^2(X')}\oplus{\rm id}_{(H^0\oplus H^4)(X')}$ provides us with an orientation preserving Hodge isometry
$\widetilde H(X,\ZZ)\congpf\widetilde H(X',\ZZ)$. Such a
composition can then be lifted to a Fourier--Mukai equivalence
$\Psi:\Db(X)\congpf\Db(X')$ (see e.g.\ \cite[Cor.\ 10.13]{HFM}).
Then $$\Psi^{-1}\circ\Phi:\Db(X)\congpf\Db(X)$$ would be an exact
equivalence with orientation reversing action on $\widetilde
H(X,\ZZ)$, which contradicts Theorem 2.
\end{proof}


\medskip

{\small\noindent {\bf Acknowledgements.} We wish to thank  M.\
Lieblich and M.\ Rapoport for useful discussions and the referees for many insightful comments and suggestions. We gratefully
acknowledge the support of the following institutions: Hausdorff
Center for Mathematics, IHES, Imperial College, Istituto Nazionale
di Alta Matematica, Max--Planck Institute, and SFB/TR 45.}



\begin{thebibliography}{99}

\bibitem{Beau} A.\ Beauville, J.-P.\ Bourguignon, M.\ Demazure,
\emph{G\'eom{\'e}trie des surfaces K3: modules et p{\'e}riodes},
S{\'e}minaires Palaiseau.  Ast{\'e}risque {\bf 126} (1985).

\bibitem{Borcea} C.\ Borcea,
\emph{Diffeomorphisms of a K3 surface}, Math.\ Ann.\ {\bf 275} (1986), 1--4.

\bibitem{B} T.\ Bridgeland,
\emph{Stability conditions on K3 surfaces}, Duke Math.\ J.\ {\bf 141} (2008), 241--291.

\bibitem{BF1} R.-O.\ Buchweitz, H.\ Flenner,
\emph{Global Hochschild (co-)homology of singular spaces}, Adv.\
Math.\ {\bf 217} (2008), 205--242.

\bibitem{BF2} R.-O.\ Buchweitz, H.\ Flenner,
\emph{The global decomposition theorem for Hochschild
(co-)homology of singular spaces via the Atiyah Chern character},
Adv.\ Math.\ {\bf 217} (2008), 243--281.

\bibitem{CvdB} D.\ Calaque, M.\ Van den Bergh,
\emph{ Hochschild cohomology and Atiyah classes},
arXiv:0708.2725.


\bibitem{Cal2} A.\ C\u{a}ld\u{a}raru,
\emph{The Mukai pairing II: The Hochschild--Kostant--Rosenberg
isomorphism}, Adv.\ Math.\ {\bf 194} (2005), 34--66.

\bibitem{Cal1} A.\ C\u{a}ld\u{a}raru, S. Willerton,
\emph{The Mukai pairing, I: a categorical approach},
arXiv:0707.2052.

\bibitem{DTT1} V. Dolgushev, D. Tamarkin, B. Tsygan, \emph{The homotopy Gerstenhaber algebra of Hochschild cochains of a regular algebra is formal}, J.\ Noncommut.\ Geom.\ {\bf 1} (2007), 1--25.

\bibitem{DTT2} V. Dolgushev, D. Tamarkin, B. Tsygan, \emph{Formality of the homotopy calculus algebra of Hochschild (co)chains}, arXiv:0807.5117.

\bibitem{Donaldson} S.\ Donaldson,
\emph{Polynomial invariants for smooth four-manifolds}, Topology
{\bf 29} (1990), 257--315.

\bibitem{Fujiki} A.\ Fujiki,
\emph{On the de Rham cohomology group of a compact K\"ahler
symplectic manifold}, Alg.\ geom., Sendai, 1985, Adv.\ Stud.\ Pure
Math.\ {\bf 10} (1987), 105--165.

\bibitem{GHJ} M.\ Gross, D.\ Huybrechts, D.\ Joyce,
\emph{Calabi--Yau manifolds and Related Geometries}, Universitext,
Springer (2003).

\bibitem{HartAG} R.\ Hartshorne,
\emph{Algebraic Geometry}, GTM 52, Springer (1977).

\bibitem{HOY} S.\ Hosono, B.H.\ Lian, K.\ Oguiso, S.-T.\ Yau,
\emph{Autoequivalences of derived category of a K3 surface and
monodromy transformations}, J.\ Alg.\ Geom.\ {\bf 13} (2004), 513--545.

\bibitem{HL} D.\ Huybrechts, M.\ Lehn,
\emph{The geometry of moduli spaces of sheaves}, Aspects of Math.\
E31, Vieweg, Braunschweig, (1997).

\bibitem{HFM} D.\ Huybrechts,
\emph{Fourier--Mukai transforms in algebraic geometry}, Oxford
Mathematical Monographs (2006).

\bibitem{HMS} D.\ Huybrechts, E.\ Macr\`i, P.\ Stellari,
\emph{Stability conditions for generic K3 surfaces}, Compositio
Math. {\bf 144} (2008), 134--162.

\bibitem{HMS1} D.\ Huybrechts, E.\ Macr\`i, P.\ Stellari,
\emph{Derived equivalences of K3 surfaces and orientation},
arXiv:0710.1645v2.

\bibitem{HMS3} D.\ Huybrechts, E.\ Macr\`i, P.\ Stellari,
\emph{Formal deformations and their categorical general fibre},
arXiv:0809.3201.

\bibitem{HS} D.\ Huybrechts, P.\ Stellari,
\emph{Equivalences of twisted K3 surfaces}, Math.\ Ann.\ {\bf 332}
(2005), 901--936.

\bibitem{HT} D.\ Huybrechts, R.\ Thomas,
\emph{Deformation-obstruction theory for complexes via Atiyah and
Kodaira--Spencer classes}, arXiv:0805.3527.

\bibitem{IllFGA} L.\ Illusie,
\emph{Grothendieck's existence theorem in formal geometry}, in
Fundamenal Algebraic Geometry, Grothendieck's FGA explained. ed.
B.\ Fantechi et al. Math.\ Surveys and Mon.\ {\bf 123}, AMS
(2005).

\bibitem{Kont} M.\ Kontsevich,
\emph{Deformation quantization  of Poisson manifolds}, Lett.\
Math.\ Phys.\ {\bf 66} (2003), 157--216.

\bibitem{Lieblich} M.\ Lieblich,
\emph{Moduli of complexes on a proper morphism}, J.\ Alg.\ Geom.\
{\bf 15} (2006), 175--206.

\bibitem{MNS} E.\ Macr\`i, M. Nieper--Wisskirchen, P. Stellari,
\emph{The module structure of Hochschild homology in some
examples}, C.\ R.\ Acad.\ Sci.\ Paris, Ser. I {\bf 346} (2008),
863-866.

\bibitem{MS} E.\ Macr\`i, P. Stellari,
\emph{Infinitesimal derived Torelli theorem for K3 surfaces}, to appear in: Int. Math. Res. Not.,
arXiv:0804.2552.

\bibitem{Markarian2} N.\ Markarian,
\emph{Poincar\'e--Birkhoff--Witt isomorphism, Hochschild homology
and Riemann--Roch theorem}. Preprint MPI 2001-52.

\bibitem{Mu} S.\ Mukai,
\emph{On the moduli space of bundles on K3 surfaces, I}, In:
Vector Bundles on Algebraic Varieties,\ Oxford University Press,\
Bombay and London (1987), 341--413.

\bibitem{Or1} D.\ Orlov,
\emph{Equivalences of derived categories and K3 surfaces}, J.\
Math.\ Sci.\ {\bf 84} (1997), 1361--1381.

\bibitem{P} D.\ Ploog,
\emph{Groups of autoequivalences of derived categories of smooth
projective varieties}, PhD-thesis, Berlin (2005).

\bibitem{Ramadoss} A.\ Ramadoss,
\emph{The relative Riemann--Roch theorem from Hochschild homology}, New York J.\ Math. {\bf 14} (2008), 643--717..

\bibitem{Sz} B.\ Szendr\H{o}i,
\emph{Diffeomorphisms and families of Fourier--Mukai transforms in
mirror symmetry}, Applications of Alg.\ Geom.\ to Coding Theory,
Phys.\ and Comp.\ NATO Science Series. Kluwer (2001), 317--337.

\bibitem{Toda} Y.\ Toda,
\emph{Deformations and Fourier--Mukai transforms}, J.\ Diff.\ Geom. {\bf 81} (2009), 197--224.

\bibitem{Weibel} C.\ Weibel,
\emph{Cyclic homology for schemes}, Proc.\
Amer.\ Math.\ Soc.\ {\bf 124} (1996), 1655--1662.

\bibitem{Yek2} A.\ Yekutieli,
\emph{The continuous Hochschild cochain complex of a scheme},
Canadian J.\ Math.\ {\bf 54} (2002), 1319--1337.

\end{thebibliography}
\end{document}